\documentclass[10pt,reqno,letter]{amsart}

\title[Analyticity and Gevrey class regularity for the Euler Equations in a bounded domain]{On the analyticity and Gevrey class regularity up to the boundary for the Euler Equations}
\author{Igor Kukavica}
\author{Vlad Vicol}

\address{Department of Mathematics, University of Southern California, Los Angeles, CA 90089}
\email{kukavica@usc.edu}
\address{Department of Mathematics, University of Southern California, Los Angeles, CA 90089}
\email{vicol@usc.edu}

\usepackage{amsmath, amsthm, amssymb}
\usepackage{color}

\usepackage[pdftex]{graphicx}
\usepackage{epstopdf}

\theoremstyle{plain}
\newtheorem{theorem}{Theorem}[section]
\newtheorem{definition}[theorem]{Definition}
\newtheorem{lemma}[theorem]{Lemma}

\theoremstyle{definition}
\newtheorem{remark}[theorem]{Remark}
\newtheorem*{ack}{Acknowledgment}

\renewcommand{\tilde}{\widetilde}

\def\GGS{{\mathcal G}^{s}}
\def\GGA{{\mathcal G}^{1}}

\def\curl{\mathop{\rm curl} \nolimits}
\def\ddiv{\mathop{\rm div} \nolimits}

\begin{document}


\begin{abstract}
We consider the Euler equations in a three-dimensional Gevrey-class bounded domain. Using Lagrangian coordinates we obtain the Gevrey-class persistence of the solution, up to the boundary, with an explicit estimate on the rate of decay of the Gevrey-class regularity radius.
\end{abstract}




\maketitle


\section{Introduction}\label{sec:intro}\setcounter{equation}{0}

The Euler equations for the velocity vector field
$u(x,t)$ and the scalar pressure field $p(x,t)$ are
given by
\begin{align}
&\partial_t u + (u\cdot \nabla) u + \nabla p = 0, \  \mbox{in}\ D \times (0,\infty), \tag{E.1}\label{eq:E1}\\
&\nabla \cdot u = 0, \  \mbox{in}\ D \times (0,\infty),  \tag{E.2}\label{eq:E2}\\
& u \cdot n = 0, \  \mbox{on}\ \partial D \times (0,\infty)
\tag{E.3}\label{eq:E3},
\end{align}
where $D$ is an open bounded Gevrey-class $s$ domain in ${\mathbb R}^3$, and $n$ is the outward unit normal to $\partial D$. We consider the initial value problem associated to
\eqref{eq:E1}--\eqref{eq:E3} with a divergence free Gevrey-class $s$ initial datum, with $s\geq 1$, namely
\begin{align}
u(0) = u_0,\tag{E.4}\  \mbox{in}\ D. \label{eq:E4}
\end{align}
The existence of smooth solutions to \eqref{eq:E1}--\eqref{eq:E4} is classical (cf.~\cite{BoB,EM,Ka,T,Y}). While in the two-dimensional case smooth initial data yield global solutions, in the three-dimensional case if $u_0 \in H^r(D)$, with $r>5/2$, the maximal time of existence of the Sobolev solution, $T_*$, might be a priori finite. If $T_* < \infty$, the vorticity must accumulate in the sense that $\int_{0}^{T_{*}} \Vert \curl u(t) \Vert_{L^\infty(D)}=\infty$ (cf.~\cite{BKM,F}). Lastly, upper bounds on  $\Vert u(\cdot,t) \Vert_{H^r(D)}$ are worse than those on $\Vert  u (\cdot,t)\Vert_{W^{1,\infty}(D)}$ due to the log-Sobolev inequality. We refer the reader to \cite{BT1,Ch,C2,MB} for the precise formulation of the above statements and for further results cornering the Euler equations.

In the present article we address the persistence of Gevrey-class regularity of the solution, i.e., we prove that if $u_0$ is of Gevrey-class $s$, then the unique Sobolev solution $u(\cdot,t) \in C([0,T_*); H^r(D))$ is of Gevrey-class $s$ for all $t< T_*$. Moreover, we are interested in sharp lower bounds on the rate of decay of the radius of analyticity and Gevrey-class regularity of the solution. We emphasize that the size of the uniform Gevrey-class radius of the solution provides an estimate for the minimal scale in dissipative flows, that is, the scale below which the Fourier coefficients decay exponentially \cite{HKR,K1}; it moreover gives the rate of this exponential decay \cite{FT,HKR}. We note that the shear flow example of Bardos and Titi~\cite{BT2} (cf.~\cite{DM}) may be used to construct explicit solutions to the three-dimensional Euler equations whose radius of analyticity, or even more generally the Gevrey-class radius, decays for all time (cf.~Remark~\ref{rem:shearflow} below).

First we summarize the rich history of this problem.
\begin{enumerate}
  \item The persistence of $C^\infty$ regularity (cf.~Foias, Frisch, and Temam~\cite{FFT}) and of real-analyticity (cf.~Bardos and Benachour~\cite{BB}) holds in both two and three dimensions.
  \item In the two-dimensional analytic case, Bardos, Benachour, and Zerner~\cite{BBZ} show that the radius of analyticity $\tau(t)$ of the solution $u(\cdot,t)$ is bounded from below as $\tau(t) \geq \exp(-C\exp(C t))/C$, for some sufficiently large constant $C$ depending on the initial data. Their elegant proof is based on analyzing the complexified equations in vorticity form.
  \item In the three-dimensional analytic case, the persistence of analyticity is proven by Bardos and Benachour~\cite{BB} using an implicit argument. In \cite{B,Be,BG1,BG2,D} using a nonlinear variant of the Cauchy-Kowalevski theorem, the authors prove the local in time existence of globally (in space) analytic solutions, with an explicit lower bound on the radius of analyticity which vanishes in finite time (independent of $T_*$). See also \cite{LCS,SC} for the dissipative Prandtl boundary layer equations.

      The proof of \cite{BB} may be modified to yield an explicit rate of decay of the radius of analyticity $\tau(t)$ which depends {\it exponentially} on $\Vert u(\cdot,t) \Vert_{H^r}$. Using different methods, Alinhac and Metivier~\cite{AM1,AM2} for the interior, and Le Bail~\cite{Lb} for the boundary value problem, obtain the short time propagation of local analyticity, with lower bounds for $\tau(t)$ that also decay {exponentially} in $\Vert u(\cdot,t) \Vert_{H^r}$. Note that the lower bounds for $\tau(t)$ obtained in \cite{AM1,AM2,BB,Lb} do not recover the lower bounds of \cite{BBZ} in the two-dimensional case, since the presently known upper bounds on high Sobolev norms of the solution increase as $C \exp(C \exp(Ct))$, for some $C>0$. Moreover, the methods used in \cite{AM1,AM2,BB,BG1,BG2,D,Lb} explicitly use the special properties of complex holomorphic functions, and hence may not be applied to the non-analytic Gevrey-class case.
  \item For the non-analytic Gevrey-class case, on a periodic domain, in both two and three dimensions, the persistence of Gevrey-class regularity follows from the elegant proof of Levermore and Oliver \cite{LO}. Their proof builds on the Fourier-based method introduced by Foias and Temam~\cite{FT} for the Navier-Stokes equations. The lower bound for the radius of Gevrey-class regularity obtained in \cite{LO} also decays {exponentially} in $\Vert u(\cdot,t) \Vert_{H^r}$. This bound was improved by the authors of the present paper in \cite{KV1}, by proving that the radius of Gevrey-class regularity decays {\it algebraically} in a high Sobolev norm of the solution, and {exponentially} in $\int_{0}^{t} \Vert \nabla u(\cdot,s) \Vert_{L^\infty}\; ds$. Therefore, the Fourier-based method may be employed (cf.~\cite{KV1}) to recover the bounds of \cite{BBZ}. For further results on analyticity cf.~\cite{Bi,BGK1,BGK2,CTV,FTi,GK1,GK2,K1,K2,KTVZ,OT}.
  \item The only result in the non-analytic Gevrey-class case, on domains with boundary, was obtained by the authors in~\cite{KV2} for $D$ a half-space. As opposed to the periodic case, here the main difficulty arises from the equation for the pressure. The classical methods of~\cite{LM,MN} are not sufficient to prove that the pressure has the same radius of Gevrey-class regularity as the velocity. In \cite{KV2} we overcome this by defining suitable norms that combinatorially encode the transfer of normal to tangential derivative in the elliptic estimate for the pressure.

      The proof of \cite{KV2} does not apply directly to the case when $D$ is a general bounded domain of Gevrey-class $s$. The main obstruction is that if $s>1$, under composition with a Gevrey-class (or even analytic) boundary straightening map, the Gevrey-class regularity radius of the velocity may deteriorate (cf.~\cite{CS,KP} and Remark~\ref{rem:FaadiBruno} below). As a consequence, we need to localize the equation using particle trajectories and define suitable Lagrangian Gevrey-class norms. This gives rise to additional difficulties because the pressure is the solution of an elliptic Neumann problem (cf.~\cite{T}), and hence is non-local.

\end{enumerate}
    The following is our main theorem.
\begin{theorem}\label{thm:main:intro}
  Let $u_0$ be divergence-free and of Gevrey-class $s$ on $D$, a Gevrey-class $s$, open bounded domain in ${\mathbb R}^3$, where $s\geq 1$, and let $r\geq 5$. Then the unique solution $u(\cdot,t) \in C([0,T_*);H^r(D))$ to the initial value problem \eqref{eq:E1}--\eqref{eq:E4} is of Gevrey-class $s$ for all $t<T_*$, where $T_* \in (0,\infty]$ is the maximal time of existence in $H^r(D)$. Moreover, the radius $\tau(t)$ of Gevrey-class regularity of the solution $u(\cdot,t)$ satisfies
\begin{align}\label{eq:thm:main:intro}
  \tau(t)\geq  C \tau_0  \exp\left(- C \left(\int_{0}^{t} \Vert u(s) \Vert_{W^{1,\infty}(D)} ds\right)^2 \right) \exp\Big(- C_0 t - C t^2 \Vert u_0 \Vert_{H^r(D)}^2\Big),
\end{align}for all $t<T_*$, where $C$ is a sufficiently large constant depending only on the domain $D$, $\tau_0$ is the radius of Gevrey-class regularity of $u_0$, and $C_0$ has additional dependence on the Gevrey-class norm of $u_0$.
\end{theorem}
\begin{remark}
   In the proof of Theorem~\ref{thm:main:intro} we also address the local (in space) propagation of Gevrey-class regularity of the solution (cf.~Theorem~\ref{thm:short:time} below), in the interior of the smooth domain, or in the neighborhood of a point where $\partial D$ is locally of Gevrey-class $s$. This extends the results of \cite{AM1,Lb} to the non-analytic Gevrey-classes.
\end{remark}
\begin{remark}\label{rem:shearflow}
Note that there exist explicit examples of solutions to \eqref{eq:E1}--\eqref{eq:E4} whose radius of Gevrey-class $s$ regularity, where $s\geq 1$, decays for all time, and vanishes as $t\rightarrow \infty$. Namely, consider the three-dimensional shear flow example (cf.~Bardos and Titi~\cite{BT2}, DiPerna and Majda~\cite{DM}) given by
\begin{align}\label{eq:Bardos:Titi}
  u(x,t) = ( f(x_2), 0 , g(x_1 - t f(x_2)) ),
\end{align}which is divergence free and satisfies \eqref{eq:E1} for smooth functions $f$ and $g$.

In the analytic category $s=1$ we may let $f(x) = \sin(x)$ and $g(x) = 1/(\tau_0^2 + \cos^2(x))$, where for simplicity $D$ is the periodic box $[0,2\pi]^3$. Substituting these particular functions $f$ and $g$ into \eqref{eq:Bardos:Titi} we obtain that $u(\cdot,0)$ has radius of analyticity $\tau_0$, while $u(x,t) = ( \sin(x_2), 0 , 1/ ( \tau_0^2 + \cos^2(x_1 - t \sin(x_2))$ has uniform radius of analyticity that decreases with the rate $1/t$. For a similar example in the non-analytic Gevrey-classes, $s>1$, let $D={\mathbb R}^3$ and define $g(x)=\exp\left(- |x|^{-1/(s-1)}\right)$ cf.~\cite{Le}. Note that $g(x)$ is of Gevrey-class $s$, but not analytic.
\end{remark}

\subsection*{Organization of the paper} In Section~\ref{sec:notation} we introduce the notation used to define the Lagrangian Gevrey-class norms. Section~\ref{sec:short-time-Gevrey} consists of a priori estimates needed to prove the short time propagation of
local analyticity (cf.~Theorem~\ref{thm:short:time}). Lemmas~\ref{lemma:velocitycommutator} and \ref{lemma:pressureterm} are proven in Sections~\ref{sec:commutator} and \ref{sec:pressure} respectively. Lastly, in Section~\ref{sec:presistence} we show how the local in space and time results may be patched together to obtain the global persistence of Gevrey-class regularity (cf.~Theorem \ref{thm:main}).

\section{Notation and preliminary remarks}\label{sec:notation}\setcounter{equation}{0}
The existence of a unique $H^r$ solution, where $r>5/2$, on a maximal time interval $[0,T_*)$, where $T_* \in (0,\infty]$, implies the existence and uniqueness of the particle trajectories (cf.~\cite{C1,MB}), that is solutions to
\begin{align}
  &\frac{d}{dt} X(t) = u (X(t),t) \label{eq:X1}\\
  &X(0) = a,\label{eq:X2}
\end{align}where $a\in \bar{D}$. For simplicity we denote by $\phi_t(a)$ the solution of \eqref{eq:X1}--\eqref{eq:X2}. It is well known that for all $t<T_*$ the maps $\phi_t \colon D \mapsto D$, and $\phi_t|_{\partial D} \colon \partial D \mapsto \partial D$ are diffeomorphisms.

\subsection*{Local change of coordinates}
Fix $x_0 \in \partial D$. In a sufficiently small neighborhood of $x_0$, the boundary $\partial D$ is the graph of a Gevrey-class $s$ function $\gamma$, i.e., for $0<r_0\ll 1$ we have $D_{r_0,x_0} = D \cap B_{r_0}(x_0) = \{ x\in B_r(x_0) \colon x_3 > \gamma(x_1,x_2)\}$. Moreover, since the Euler equations are invariant under rigid body rotations of ${\mathbb R}^3$, modulo composition with a rigid body rotation about $x_0$, we may assume that $\Vert \partial_1 \gamma \Vert_{L^\infty(\bar{D}_{r_0,x_0}')} + \Vert \partial_2 \gamma \Vert_{L^\infty(\bar{D}_{r_0,x_0}')} \leq \overline{\varepsilon} \ll 1$, for $r_0$ sufficiently small, where $\overline{\varepsilon}$ is a fixed, sufficiently small universal constant, to be chosen later. Here we have denoted $D_{r_0,x_0}' = \{ x' \colon x \in D_{r_0,x_0}\}$, where we write $x'=(x_1,x_2)$ for $x=(x_1,x_2,x_3)$. Define a boundary straightening map $\theta: {\mathbb R}^3 \rightarrow {\mathbb R}^3$ by
\begin{align}
\theta(x_1,x_2,x_3) = (x_1,x_2,x_3 - \gamma(x_1,x_2))= (y_1,y_2,y_3).
 \end{align}Note that $\det (\partial \theta/ \partial x) = 1$. By the construction of $\theta$ we have $\tilde{D}_{r_0,x_0} = \theta(D_{r_0,x_0}) = \{ y \in \theta(B_{r_0}(x_0)) \colon y_3 >0\}$.

Let $\Omega = D\cap B_{r_0/2}(x_0) \subset D_{r_0,x_0}$ be a neighborhood of $x_0$. Also let $\tilde{\Omega} = \theta(\Omega)$ and $\tilde{\Omega}_t  = \theta(\Omega_t)$.

\begin{figure}[!h]
  \includegraphics[width=84ex]{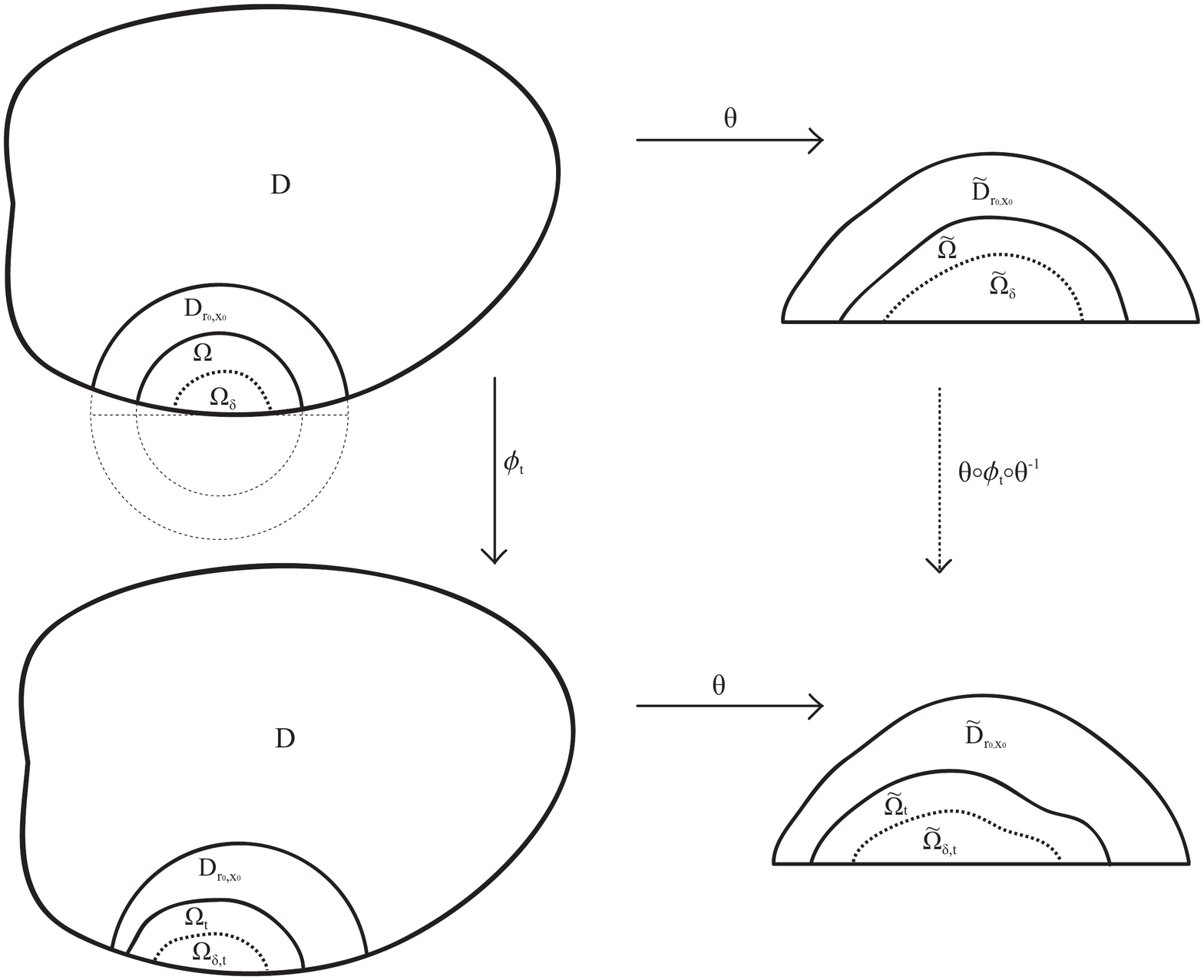}
\end{figure}

There exists $T_1 = T_1(r_0, u )$ such that for all $0=T_0 \leq t \leq T_1$ we have $\Omega_t = \phi_t (\Omega) \subset D_{r_0,x_0}$. The value of $T_1$ may be estimated from below by using the representation formula for solutions of \eqref{eq:X1}--\eqref{eq:X2}. We have
\begin{align}
| \phi_t (a) -a | \leq \int_{0}^{t} | u(\phi_s(a),s)|\; ds \leq \int_{0}^{t} \Vert u(\cdot,s) \Vert_{L^\infty(D_{r_0,x_0})}\; ds \leq K(t),
\end{align}where we set
\begin{align}
  K(t) = \int_{0}^{t} \Vert u(s) \Vert_{W^{1,\infty}(D)}\; ds. \label{eq:K:def}
\end{align}Therefore, it is sufficient to chose $T_1$ such that $K(T_1) \leq {\rm {dist}}(\bar{\Omega},\partial B_{r_0}(x_0)) = r_0/2$.

On the closure of $\tilde{D}_{r_0,x_0}$ we let $\varrho$ be the Euclidean distance to the curved part of the boundary of $\tilde{\Omega}$, that is, $\varrho(y) = 0$ if $y \in \tilde{\Omega}^c$ and $\varrho(y) = {\rm {dist}}(y,\partial \tilde{\Omega} \setminus \{ y_3 = 0\})$ if $y \in \tilde{\Omega}$. As in \cite{AM1,AM2,Lb}, for all $0<\delta \leq \delta_0$, where $\delta_0 > 0$ is sufficiently small, we define the set
\begin{align}
\tilde{\Omega}_\delta = \{ y \in \tilde{\Omega}\colon \varrho(y) > \delta\}.
 \end{align}By the triangle inequality and the definition of $\varrho$ it follows that $|y^{(1)}-y^{(2)}|\geq r$ for all $y^{(1)}\in \tilde{\Omega}_{\delta + r}$ and $y^{(2)}\in \Omega_{\delta}^c$. Also let $\Omega_\delta = \theta^{-1}(\tilde{\Omega}_\delta)$, $\Omega_{t,\delta} = \phi_t(\Omega_\delta)$ and $\tilde{\Omega}_{t,\delta} = \theta(\Omega_{t,\delta})$. Here $\delta_0 = \delta_0(\gamma) \leq 1$ is chosen small enough such that for all  $\delta \in [0,\delta_0)$, the set $\Omega_\delta$ is a Gevrey-class $s$ domain, i.e., it lies on one side of a Gevrey-class surface.


If $y^{(1)}\in \tilde{\Omega}_{\delta+r,t}$ and $y^{(2)} \in \tilde{\Omega}_{\delta,t}$, where $\delta, r+\delta \in (0,\delta_0)$, it follows by the mean value theorem that
\begin{align}
  r \leq |\theta\circ \phi_{-t}\circ \theta^{-1} (y^{(1)}) - \theta\circ \phi_{-t}\circ \theta^{-1}(y^{(2)}) | &\leq C |y^{(1)}-y^{(2)}| \Vert \nabla \phi_{-t} \Vert_{L^\infty(D)}\notag\\
  &\leq  C |y^{(1)}-y^{(2)}| (1 + K^2(t)),\label{eq:inversematrix}
\end{align} where $C$ is a constant depending on $\theta$. In \eqref{eq:inversematrix} we have used that $\nabla \phi_{-t}$ is the inverse matrix of $\nabla \phi_t$ (whose determinant is $1$ since $\ddiv u =0$), and the fact that the $2\times 2$ minors of this matrix are bounded by $1+K^2$. Therefore, by \eqref{eq:inversematrix}, we have $|y^{(1)}-y^{(2)}| \geq r/ (C + C K^2(t))$. Hence there exists a smooth cut-off function $\eta$ such that $\eta\equiv 1$ on $\tilde{\Omega}_{\delta+r,t}$ and $\eta \equiv 0$ on $\tilde{\Omega}_{\delta,t}^c$, with
\begin{align}
  |\nabla \eta| \leq \frac{C+ C K^{2}(t)}{r}, \label{eq:cutoff:estimate}
\end{align}for some positive constant $C=C(D)$. We denote $\tilde{u}(y,t) = u(x,t)$ and similarly $\tilde{p}(y,t) = p(x,t)$.


\subsection*{Gevrey-class norms}
We recall (cf.~\cite{KP,Le,LM}) the definition of the Gevrey-class $s$, denoted by $\GGS$.
\begin{definition}\label{def:GGS}
A function $v\in C^\infty(D)$ is said to be of Gevrey-class $s$ on $D$, where $s\geq 1$, written $v \in \GGS$, if there exist positive constants $M,\tau>0$ such that
\begin{align}\label{eq:GGS:def}
\Vert \partial^\alpha v \Vert_{L^\infty(D)} \leq M \frac{|\alpha|!^s}{\tau^{|\alpha|}}
\end{align}for all multi-indices $\alpha \in {\mathbb N}_{0}^{3}$. We refer to the constant $\tau$ in \eqref{eq:GGS:def} as the radius of Gevrey-class regularity of $v$, or simply as the $\GGS$-radius of $v$.
\end{definition} As opposed to the class of real analytic functions $\GGA$, functions in $\GGS$ with $s>1$ may have compact support, they may vanish of infinite order at a point, and there exist $\GGS$ partitions of unity (cf.~\cite{KP}). The $\GGS$ norms used in this paper are defined as follows. For a Gevrey-class $s$ function $\tilde{v}(y,t)$ denote
\begin{align}
  [\tilde{v}(t)]_m =  \sum\limits_{|\alpha|=m} \epsilon^{\alpha_3} \sup\limits_{0<\delta\leq \delta_0} \delta^{m-3} \Vert \partial^\alpha \tilde{v}(\cdot,t) \Vert_{L^2(\tilde{\Omega}_{t,\delta})},\label{eq:[]m}
\end{align}for all $m\geq 3$. In this paper we work with the Lagrangian Gevrey-class $s$ norm defined by
\begin{align}
  \Vert \tilde{v}(t) \Vert_{X_{\tau(t)}} &= \sum\limits_{m=3}^{\infty} [\tilde{v}(t)]_m \frac{\tau(t)^{m-3}}{(m-3)!^s}\label{eq:Xtaudef},
  \end{align}where $s\geq 1$, and $\tau>0$. We also let
  \begin{align}
      \Vert \tilde{v}(t) \Vert_{Y_{\tau(t)}} &= \sum\limits_{m=4}^{\infty} [\tilde{v}(t)]_m \frac{m \tau(t)^{m-4}}{(m-3)!^s}\label{eq:Ytaudef}.
  \end{align}
\begin{remark}\label{rem:FaadiBruno}
  If $\Vert \tilde{u}(y)\Vert_{X_\tau} < \infty$, it follows from the Sobolev inequality that $\tilde{u}\in\GGS$ and that $\tilde{u}(y)$ has Gevrey-class regularity radius at least $\tau$. As opposed to the analytic case, if $s>1$, the map $\theta^{-1}\colon y \mapsto x$ possibly shrinks the radius by a constant factor $0< a_*\leq 1$, where $a_* = a_*(\gamma)$. This fact may be proven using the multi-dimensional generalization of the Fa\'a~di~Bruno formula (cf.~\cite{CS,KP}). Thus, if $\tilde{u}(y)$ has Gevrey-class radius $\tau$, then $u(x)$ has $\GGS$-radius at least $a_* \tau$.
\end{remark}
\subsection*{Notation} When it is clear from the context that we are working with a function on the flattened domain, we simply write $v$ instead of $\tilde{v}$. In the present paper we set $n! = 1$ whenever $n\leq 0$. Also we use the notation $\Vert D^k v \Vert_{L^p} = \sum_{|\alpha|=k} \Vert  \partial^\alpha v \Vert_{L^p}$ and similarly $\Vert {D'^k} v \Vert_{L^p} = \sum_{|\alpha|=k, \alpha_3=0} \Vert  \partial^\alpha v \Vert_{L^p}$. Lastly, $C$ denotes a sufficiently large positive constant which may depend on the domain.

\section{Short time local Gevrey-class a priori estimates} \label{sec:short-time-Gevrey}\setcounter{equation}{0}
The proof of Theorem~\ref{thm:main:intro} consists of a priori estimates. These estimates can be made rigorous by noting that $u(\cdot,t)\in C^\infty(\bar{D})$ for all $t<T_*$ (cf.~\cite{FFT}), and by performing all below estimates on truncated sums $\sum_{m=3}^{q} [\tilde{v}]_m \tau^{m-3}/(m-3)!^s$. For $q\geq 5$ these estimates close, since the energy estimates for the Euler equations close in Sobolev spaces (cf.~\cite{T}), and are independent of $q$, so we may let $q\rightarrow \infty$. 

Let $d^+ f(t)/dt = \limsup_{h\rightarrow 0+} (f(t+h)-f(t))/h$ denote the right derivative of a function $f(t)$, which agrees with the usual derivative if the latter exists. Using the definitions \eqref{eq:[]m}--\eqref{eq:Xtaudef} we obtain
\begin{align}
  \frac{d^+}{dt} \Vert \tilde{u}(t) \Vert_{X_{\tau(t)}} \leq \dot{\tau}(t) \Vert \tilde{u}(t) \Vert_{Y_{\tau(t)}} + \sum\limits_{m=3}^{\infty} \left( \sum\limits_{|\alpha|=m}\epsilon^{\alpha_3} \frac{d^+}{dt} \sup\limits_{0<\delta\leq \delta_0} \delta^{m-3}  \Vert \partial^{\alpha}  \tilde{u}(t) \Vert_{L^2(\tilde{\Omega}_{\delta,t})}\right) \frac{\tau(t)^{m-3}}{(m-3)!^s}. \label{eq:ode1}
\end{align}In order to switch the $d^+/dt$ and the $\sup_\delta$ (cf.~Lemma~\ref{lemma:ap:3}) we need upper bounds for $(d/dt) \Vert \partial^\alpha \tilde{u}(t) \Vert_{\tilde{\Omega}_{\delta,t}}$ for all $|\alpha|\geq 3$. The following lemma is a Lagrangian energy estimate in the straightened domain and provides the desired upper bound.
\begin{lemma}\label{lemma:energyestimate}
  For all $\alpha \in {\mathbb N}_{0}^{3}, t>0$, and $0<\delta\leq \delta_0$, we have
  \begin{align}\label{eq:energyestimate}
  \frac{d}{dt} \Vert \partial^\alpha \tilde{u}(\cdot,t) \Vert_{L^2(\tilde{\Omega}_{\delta,t})} \leq \Vert [\partial^\alpha,\tilde{u}_j\, \partial_j \theta_k\, \partial_k]\tilde{u}(\cdot,t) \Vert_{L^2(\tilde{\Omega}_{\delta,t})} + \Vert \partial^\alpha (\partial_j \theta_k\, \partial_k \tilde{p}(\cdot,t)) \Vert_{L^2(\tilde{\Omega}_{\delta,t})} = M_\delta(t) ,
  \end{align}where the bracket $[\cdot,\cdot]$ represents a commutator.
\end{lemma}
\begin{proof}The standard Lagrangian energy estimate (cf.~\cite[Lemma 2.3]{AM1} and \cite[Section 2.b]{Lb}) shows that a smooth solution $v$ to $\partial_t v + (u\cdot \nabla) v = g$ satisfies
\begin{align}
  \frac{d}{dt} \Vert v(t,\cdot) \Vert_{L^2(\Omega_{\delta,t})}^2 = \frac{d}{dt} \int\limits_{\Omega_\delta} |v(t,\phi_t(x))|^2\, dx = 2 \int\limits_{\Omega_{\delta,t}} g(t,x) v(t,x)\, dx
\end{align}Here we used the fact that $\ddiv u =0$ implies $\det( \partial \phi_t(x) / \partial x) = 1$. Let $y = \theta(x)$ and denote $\tilde{v}(y) = v(x)$. Similarly define $\tilde{u}(y) = u(x)$ and $\tilde{g}(y) = g(x)$. Then $\tilde{v}$ solves the equation $\partial_t \tilde{v} + \tilde{u}_j\; \partial_j \theta_k\; \partial_k \tilde{v} = \tilde{g}$, and since $\det(\partial \theta/\partial x) = 1$, we have $\Vert\tilde{v}(t,\cdot)\Vert_{L^2(\tilde{\Omega}_{\delta,t})} = \Vert v(t,\cdot) \Vert_{L^2(\Omega_{\delta,t})}$. The lemma follows from the above remarks with $\tilde{v}= \partial^\alpha \tilde{u}$ and $\tilde{g} = [\partial^\alpha, \tilde{u}_j\, \partial_j \theta_k\,  \partial_k] \tilde{u} + \partial^\alpha( \partial_j \theta_k\, \partial_k \tilde{p})$, and the H\"older inequality.
\end{proof}
Using the bound \eqref{eq:energyestimate}, from Lemma~\ref{lemma:ap:3} we obtain
\begin{align*}\frac{d^+}{dt} \sup_{0<\delta\leq\delta_0} \delta^{m-3} \Vert \partial^\alpha \tilde{u}(t)\Vert_{\tilde{\Omega}_{\delta,t}} \leq \sup_{0<\delta\leq \delta_0} \delta^{m-3} M_\delta(t).\end{align*}Therefore, we may estimate the sum on the right of \eqref{eq:ode1} as
\begin{align}
  \frac{d^+}{dt} \Vert \tilde{u}(t) \Vert_{X_{\tau(t)}} &\leq \dot{\tau}(t) \Vert \tilde{u}(t) \Vert_{Y_{\tau(t)}} + \sum\limits_{m=3}^{\infty} \left( \sum\limits_{|\alpha|=m}\epsilon^{\alpha_3} \sup\limits_{0<\delta\leq \delta_0} \delta^{m-3}  M_\delta(t) \right) \frac{\tau(t)^{m-3}}{(m-3)!^s}\notag\\
  & \leq \dot{\tau}(t) \Vert \tilde{u}(t) \Vert_{Y_{\tau(t)}} + {\mathcal C} + {\mathcal P} \label{eq:ode2},
\end{align}where
\begin{align}
  {\mathcal C} = \sum\limits_{m=3}^{\infty} \left( \sum\limits_{|\alpha|=m}\epsilon^{\alpha_3} \sup\limits_{0<\delta\leq \delta_0} \delta^{m-3}  \Vert [\partial^\alpha, \tilde{u}_j\, \partial_j \theta_k\,  \partial_k]\tilde{u}(\cdot,t) \Vert_{L^2(\tilde{\Omega}_{\delta,t})}\right) \frac{\tau(t)^{m-3}}{(m-3)!^s},\label{eq:Cdef}
\end{align}
and
\begin{align}
   {\mathcal P} = \sum\limits_{m=3}^{\infty} \left( \sum\limits_{|\alpha|=m}\epsilon^{\alpha_3} \sup\limits_{0<\delta\leq \delta_0} \delta^{m-3}  \Vert\partial^\alpha( \partial_j \theta_k\, \partial_k \tilde{p}(\cdot,t) ) \Vert_{L^2(\tilde{\Omega}_{\delta,t})}\right) \frac{\tau(t)^{m-3}}{(m-3)!^s}.\label{eq:Cdef}
\end{align}
The estimates for ${\mathcal C}$ and ${\mathcal P}$ are given in the following two lemmas.
\begin{lemma}\label{lemma:velocitycommutator} If $\tau < \tau_*$, where $\tau_*$ is the Gevrey-class regularity radius of the boundary, then the following estimate holds
\begin{align}\label{eq:velocitycommutator}
{\mathcal C} &\leq C (1+\tau^2) \Big(\Vert \tilde{u} \Vert_{W^{2,\infty}(\tilde{\Omega}_t)}^2 + \Vert \tilde{u} \Vert_{H^5(\tilde{\Omega}_t)}^2\Big)\notag\\
 &\ + C \Bigg( \tau \Vert D \tilde{u} \Vert_{L^\infty(\tilde{\Omega}_t)}  + (\tau^2+\tau^3) \Big( \Vert \tilde{u} \Vert_{W^{2,\infty}(\tilde{\Omega}_t)} + \Vert \tilde{u} \Vert_{H^5(\tilde{\Omega}_t)} \Big) + (\tau^{3/2}+ (1+K^3)\tau^3) \Vert \tilde{u} \Vert_{X_\tau} \Bigg) \Vert \tilde{u} \Vert_{Y_\tau},
\end{align}where $C$ is a sufficiently large positive constant depending on $\gamma$, and $K$ is as defined in \eqref{eq:K:def}.
\end{lemma}
The proof of Lemma~\ref{lemma:velocitycommutator} is given in Section~\ref{sec:commutator}, while the proof of Lemma~\ref{lemma:pressureterm} below is given in Section~\ref{sec:pressure}.
\begin{lemma}\label{lemma:pressureterm}
For $\epsilon >0$ fixed, sufficiently small depending only on $\gamma$, if $\tau \leq \epsilon \tau_*$, where $\tau_*$ is the Gevrey-class regularity radius of the boundary, then we have
\begin{align}\label{eq:pressureterm}
 {\mathcal P} &\leq C  (1+\tau^2) \Big( \Vert \tilde{u} \Vert_{W^{2,\infty}(\tilde{\Omega}_t)}^2 + \Vert \tilde{u}\Vert_{H^3(\tilde{\Omega}_t)}^2 + (1+K^2)\Vert \tilde{p} \Vert_{H^4(\tilde{\Omega}_t)} + \Vert \tilde{p} \Vert_{W^{3,\infty}(\tilde{\Omega}_t)}\Big)\notag\\
 &\ + C \Bigg( \tau \Vert \tilde{u} \Vert_{W^{1,\infty}(\tilde{\Omega}_t)} + (\tau^2+\tau^3) \Vert \tilde{u} \Vert_{W^{2,\infty}(\tilde{\Omega}_t)} + \tau^2 \Vert \tilde{u} \Vert_{H^5(\tilde{\Omega}_t)} + (\tau^{3/2} + (1+K^3) \tau^{4}) \Vert \tilde{u} \Vert_{X_\tau} \Bigg) \Vert \tilde{u} \Vert_{Y_\tau},
\end{align}for some sufficiently large constant $C$ depending on $\gamma$, where $K$ is as defined in \eqref{eq:K:def}.
\end{lemma}
By combining estimates \eqref{eq:ode2}, \eqref{eq:velocitycommutator}, and \eqref{eq:pressureterm}, with the Sobolev embedding, and the classical pressure estimate in Sobolev spaces $\Vert p \Vert_{H^m(D)} \leq C \Vert u \Vert_{H^{m-1}(D)}^2$ (cf.~\cite[Lemma 1.2]{T}), we obtain for $r\geq 5$ that
\begin{align}\label{eq:ode3}
\frac{d}{dt} \Vert \tilde{u} \Vert_{X_{\tau}} &\leq C (1+\tau^2)(1+K^2) \Vert u \Vert_{H^r(D)}^2 \notag\\
& + \Big( \dot{\tau} + C \tau \Vert u \Vert_{W^{1,\infty}(D)}  + C (\tau^2 +\tau^3) \Vert u \Vert_{H^r(D)}+ C (\tau^{3/2}+(1+K^3)\tau^4) \Vert \tilde{u} \Vert_{X_\tau}\Big) \Vert \tilde{u} \Vert_{Y_{\tau}}.
\end{align}for some fixed positive constant $C$ depending on the domain $D$. Let
\begin{align}
    {M}(t) = \Vert u(\cdot,t) \Vert_{H^r(D)}
\end{align}
and
\begin{align}
 {N}(t) = \Vert u(\cdot,t) \Vert_{W^{1,\infty}(D)}  \label{eq:N:def}
\end{align}for all $0\leq t < T_*$. Note that $K(t)= \int_{0}^{t} N(s)\, ds$. By possibly increasing the constant $C=C(D)$, we have
\begin{align}\label{eq:ode4}
\frac{d}{dt} \Vert \tilde{u} \Vert_{X_{\tau}} &\leq C (1+\tau^2)(1+K^2)  {M}^2 + \Big( \dot{\tau} + C \tau N  + C (\tau^{2} +\tau^3) M + C (\tau^{3/2}+(1+K^3) \tau^4) \Vert \tilde{u} \Vert_{X_\tau}\Big) \Vert \tilde{u} \Vert_{Y_{\tau}}.
\end{align}Let $\tau(t)$ be chosen such that $\tau(t) \leq \tau_0 \leq \tau_*$, where $\tau_*$ is the radius of Gevrey-class regularity of the boundary, and for all $0=T_0\leq t \leq T_1$ let $\tau(t)$ be the solution of
\begin{align}
&\dot{\tau} + C_0 \tau N  + C_0 \tau^{3/2} L = 0,\label{eq:tau:construct}
\end{align}with the initial condition $\tau(0)=\tau_0$, where $C_0$ is a sufficiently large  fixed positive constant (for instance $C_0/(2+\tau_*^{2}) > C$, the constant of \eqref{eq:ode4}); we have denoted
\begin{align}
  L(t) = C_0 M(t) + \left(1 + C_0\Big(1+K^3(t)\Big)\right) \left(\Vert \tilde{u}_0 \Vert_{X_{\tau_0}} + C_0 \int_{0}^{t}\Big(1+K^2(s)\Big) M^2(s)\; ds\right). \label{eq:L:def}
\end{align} Then $\tau$ is decreasing, and by \eqref{eq:ode4} for short time we have
\begin{align}\label{eq:Gevrey:growth}
  \Vert \tilde{u}(t) \Vert_{X_{\tau(t)}}\leq \Vert \tilde{u}_0 \Vert_{X_{\tau_0}} + C_0\int_{0}^{t} \Big(1+K^2(s)\Big)M^2(s)\; ds.
\end{align}By \eqref{eq:ode4}, if \eqref{eq:tau:construct} holds for all $t\in[T_0,T_1]$, then $\tilde{u}(t) \in X_{\tau(t)}$ and \eqref{eq:Gevrey:growth} is also valid for all $t\in [T_0,T_1]$. The radius of Gevrey-class regularity $\tau(t)$ may be computed explicitly from \eqref{eq:tau:construct} as
\begin{align}
  &\tau(t) = \exp\Bigl( - C_0 K(t) \Bigr) \Bigg( \tau_{0}^{-1/2} + C_0 \int_{0}^{t} L(s) \exp\Bigl( - C_0 K(s)\Bigr) \; ds\Bigg)^{-2}\label{eq:tau:nasty},
\end{align} where $L$ is defined by \eqref{eq:L:def}. By further estimating the Sobolev norms in \eqref{eq:tau:nasty} using
\begin{align*}
  M^2(t) = \Vert u(t) \Vert_{H^r(D)}^2 \leq C \Vert u_0 \Vert_{H^r(D)}^2 \exp\left(C \int_{0}^{t} \Vert u(s) \Vert_{W^{1,\infty}(D)}ds\right) = C \Vert u_0 \Vert_{H^r(D)}^2 e^{C K(t)},
\end{align*}for some positive constant $C=C(r,D)$, we obtain a more compact lower bound for $\tau(t)$ given by
\begin{align}
  \tau(t) &\geq \tau_0 \left( 1 + C t  \Vert \tilde{u}_0 \Vert_{X_{\tau_0}} + C t^2 \Vert u_0 \Vert_{H^r}^2 \Big(1+K^5(t)\Big)\right)^{-2} \exp\Bigl(- C K(t)\Bigr)\notag\\
  & \geq \tau_0 \left( 1 + C t  \Vert \tilde{u}_0 \Vert_{X_{\tau_0}} + C t^2 \Vert u_0 \Vert_{H^r}^2\right)^{-2} \exp\Bigl(- C K(t)\Bigr);\label{eq:tau:compact}
\end{align}we used  $(1+x^5)^{-2} \geq \exp(-2x)$ for all $x\geq 0$, where $C=C(D,r)$ is a sufficiently large positive constant. Therefore we have proven the following theorem.

\begin{theorem}\label{thm:short:time}
Let $u_0$ be divergence-free and of Gevrey-class $s$, with $s\geq 1$, on a Gevrey-class $s$, open, bounded domain $D \subset {\mathbb R}^3$. Fix $r\geq 5$, $x_0 \in \partial D$, and $r_0>0$ sufficiently small. Let $\Omega$  be a neighborhood of $x_0$ compactly embedded in $ B_{r_0}(x_0) \cap D$, and let $T_1$ be the maximal time such that $\phi_t(\Omega) \subset B_{r_0}(x_0) \cap D$ for all $0 \leq t < T_1$. Then the unique $H^r$-solution $u(\phi_t(\cdot),t)$ to the initial value problem \eqref{eq:E1}--\eqref{eq:E4} is of Gevrey-class $s$ for all $t<T_1$. Moreover, there exist constants $\epsilon= \epsilon(D)$, and $\tau_* = \tau_*(D)$, such that if $\tilde{u}(0) \in X_{\tau_0}$, and $\tau_0\leq \epsilon \tau_*$, then $\tilde{u}(\cdot,t) \in X_{\tau(t)}$ for all $t\in [0,T_1)$, where the Gevrey-class radius $\tau(t)$ of the solution $u(\phi_t(\cdot),t)$ satisfies
\begin{align}\label{eq:thm:mainbound}
  \tau(t)\geq  \tau_0 \left( 1 +Ct \Vert \tilde{u}_0 \Vert_{X_{\tau_0}} + C t^2 \Vert u_0 \Vert_{H^r}^2 \right)^{-2} \exp\left(- C \int_{0}^{t} \Vert u(s) \Vert_{W^{1,\infty}} ds\right)
\end{align}for all $t<T_1$, with $C$ a sufficiently large constant depending only on the domain $D$.
\end{theorem}
\begin{remark}\label{rem:interior}
  Theorem~\ref{thm:short:time} gives the local in time Gevrey-class persistence at the boundary of $D$. The short-time Gevrey-class persistence in the interior of $D$, with explicit bound on the radius of Gevrey-class regularity is obtained using similar arguments to the ones given in this section. Namely, given $x_0 \in D$ and $r_0>0$ sufficiently small, we let $\Omega$ be a Gevrey-class neighborhood of $x_0$, with $\Omega \subset B_{r_0}(x_0)\cap D$. The main step is to show that for all $t>0$ such that $\phi_t (\Omega) \subset B_{r_0}(x_0)\cap D$, the analogous estimates to the ones given in Lemmas~\ref{lemma:velocitycommutator} and \ref{lemma:pressureterm} hold. The bound on the velocity commutator ${\mathcal C}$ is proven by repeating exactly the same estimates as in Section~\ref{sec:commutator} below. Since we are away from the boundary, the bound for the pressure term ${\mathcal P}$ is obtained from classical interior elliptic estimates for the pressure (cf.~\eqref{eq:P.1}) and arguments parallel to the ones presented in Lemma~\ref{lemma:P:u}. Since the interior pressure estimates are only simpler than the boundary case, we omit further details. It follows that the solution $u(\cdot,t)$ is of Gevrey-class $s$ on $\phi_t(\Omega)$, the radius of Gevrey-class regularity $\tau(t)$ satisfies the lower bound \eqref{eq:thm:mainbound}, and that the Gevrey-class norm is bounded from below by the right side of \eqref{eq:Gevrey:growth}.
\end{remark}

\section{The velocity commutator estimate}\label{sec:commutator}\setcounter{equation}{0}
Since in this section we work only for a fixed time $t$ and on the straightened domain, we suppress the time dependence and the tilde for all functions and domains. The goal of this section is to prove Lemma~\ref{lemma:velocitycommutator}, that is to estimate
\begin{align*}
  {\mathcal C} &= \sum\limits_{m=3}^{\infty} \sum\limits_{|\alpha|=m} \epsilon^{\alpha_3}  \sup\limits_{0< \delta \leq \delta_0} \Big( \delta^{m-3} \Vert [\partial^{\alpha}, u_j\, \partial_j \theta_k\, \partial_k] u \Vert_{L^2(\Omega_{\delta})}\Big) \frac{\tau^{m-3}}{(m-3)!^s}.
\end{align*}
\begin{proof}[Proof of Lemma~\ref{lemma:velocitycommutator}]The proof consists of two parts. First we estimate the $\partial_j \theta_k$ coefficients from the definition of $\mathcal{C}$ and exploit the commutator (cf.~\eqref{eq:eliminate:theta:5} below). Then we estimate the Gevrey-class norm of $u_i\; \partial_j u_k$ (cf.~\eqref{eq:C1:bound}--\eqref{eq:Clow+Chigh} below) for all $1\leq i,j,k\leq 3$.

The Leibniz rule and the fact  $\sup_\delta \sum_n x_{n,\delta} \leq \sum_n \sup_\delta x_{n,\delta}$ for all sequences $x_{n,\delta}$, imply
\begin{align}
  {\mathcal C} & \leq \sum\limits_{m=3}^{\infty} \sum\limits_{|\alpha|=m} \sum\limits_{0<\beta \leq \alpha} \sum\limits_{0\leq \gamma \leq \beta} {\alpha \choose \beta} {\beta \choose \gamma} \epsilon^{\alpha_3} \Vert \partial^\gamma \partial_j \theta_k \Vert_{L^\infty(\Omega)} \sup\limits_{0<\delta \leq \delta_0} \Big( \delta^{m-3}  \Vert \partial^{\beta-\gamma} u_j\, \partial^{\alpha - \beta} \partial_k u \Vert_{L^2(\Omega_{\delta})}\Big) \frac{\tau^{m-3}}{(m-3)!^s},\label{eq:eliminate:theta:1}
\end{align}where $\Omega = \bigcup_{0<\delta\leq \delta_0} \Omega_\delta$. Since the boundary is of Gevrey-class $s$, there exist constants $C,\tau_*>0$ such that
\begin{align}
  \sum\limits_{|\beta|=n} \Vert \partial^\beta D \theta_k \Vert_{L^\infty(\Omega)} \leq C \frac{(n-3)!^s}{\tau_{*}^n},\label{eq:eliminate:theta:2}
\end{align}for all $n\geq 0$. Using ${\alpha \choose \beta} {\beta \choose \gamma} = {\alpha \choose \gamma} {\alpha-\gamma \choose \beta-\gamma} \leq {\alpha \choose \gamma} {m-k \choose j-k}$, we may rewrite the right side of \eqref{eq:eliminate:theta:1} as
\begin{align}
  {\mathcal C} &\leq C \sum\limits_{m=3}^{\infty} \sum\limits_{j=1}^{m} \sum\limits_{k=0}^{j} {m \choose k} \frac{(k-3)!^s\, (m-k-3)!^s}{(m-3)!^s} \left(\frac{\tau}{\tau_*}\right)^k\notag\\
  & \qquad \times \sum\limits_{|\alpha|=m}\sum\limits_{|\beta|=j,\, \beta\leq \alpha} \sum\limits_{|\gamma|=k,\, \gamma\leq \beta}\left({\alpha \choose \gamma}{m \choose k}^{-1} \Vert \partial^\gamma D\theta \Vert_{L^\infty(\Omega)} \frac{\tau_{*}^{k}}{(k-3)!^s} \right)\notag\\
  &  \qquad \qquad \times \left(\epsilon^{\alpha_3 } \frac{\tau^{m-k-3}}{(m-k-3)!^s} {m-k \choose j-k} \sup\limits_{0<\delta\leq \delta_0}  \delta^{m-3} \Vert \partial^{\beta-\gamma} u\Vert_{L^p(\Omega_\delta)} \Vert \partial^{\alpha-\beta} Du\Vert_{L^{2p/(p-2)}(\Omega_\delta)} \right),\label{eq:eliminate:theta:3}
  \end{align}where $p=2$ if $j-k > m-j$, and $p=\infty$ if $j-k \leq m-j$. Observe that $\tau/\tau_* < 1$. Since $s\geq 1$, there exists a constant $C$ such that
  \begin{align}\label{eq:eliminate:theta:combinatorial}
    {m \choose k}
    \frac{(k-3)!^s (m-k-3)!^s }{(m-3)!^{s}} \leq \frac{C}{(m-k+1)^{s-1}} + C \chi_{\{k=0\}}
  \end{align}for all $0\leq k \leq m$, where $\chi_{\{k=0\}} = 1$ if $k=0$, and $\chi_{\{k=0\}} = 0$ if $k\geq 1$. By \eqref{eq:eliminate:theta:2}, \eqref{eq:eliminate:theta:combinatorial}, Lemma~\ref{lemma:combinatorial}, and using ${\alpha \choose \gamma} \leq {m \choose k}$, we obtain
 \begin{align}
  {\mathcal C} &\leq C \sum\limits_{m=3}^{\infty} \sum\limits_{j=1}^{m} \sum\limits_{k=0}^{j} \left(\frac{\tau}{\tau_*}\right)^k  \sum\limits_{|\alpha|=m}\sum\limits_{|\beta|=j,\, \beta\leq \alpha} \sum\limits_{|\gamma|=k,\, \gamma\leq \beta} \left(\Vert \partial^\gamma D\theta \Vert_{L^\infty(\Omega)} \frac{\tau_{*}^{k}}{(k-3)!^s}\right)\frac{\tau^{m-k-3}}{(m-k-3)!^s} {m-k \choose j-k} \notag\\
  & \qquad \qquad  \times \left(\frac{1}{(m-k+1)^{s-1}} + \chi_{\{k=0\}}\right) \left(\epsilon^{\alpha_3 } \sup\limits_{0<\delta\leq \delta_0}  \delta^{m-3} \Vert \partial^{\beta-\gamma} u\Vert_{L^p(\Omega_\delta)} \Vert \partial^{\alpha-\beta} Du\Vert_{L^{2p/(p-2)}(\Omega_\delta)}\right)\notag\\
  & \leq C \sum\limits_{m=3}^{\infty} \sum\limits_{j=1}^{m} \sum\limits_{k=0}^{j} \left(\frac{\tau}{\tau_*}\right)^k  \sum\limits_{|\alpha|=m-k}\sum\limits_{|\beta|=j-k,\, \beta\leq \alpha}\frac{\tau^{m-k-3}}{(m-k-3)!^s} {m-k \choose j-k} \notag\\
  &  \qquad \qquad  \times \left(\frac{1}{(m-k+1)^{s-1}} + \chi_{\{k=0\}}\right) \left(\epsilon^{\alpha_3 }  \sup\limits_{0<\delta\leq \delta_0}\delta^{m-3} \Vert \partial^{\beta} u\Vert_{L^p(\Omega_\delta)}   \Vert\partial^{\alpha-\beta} Du\Vert_{L^{2p/(p-2)}(\Omega_\delta)} \right).\label{eq:eliminate:theta:4}
  \end{align}Due to the definition of the Gevrey-class norm, in \eqref{eq:eliminate:theta:4} we need to consider the cases $m-k<3$ and $m-k\geq 3$ separately. We estimate the discrete convolution using Lemma~\ref{lemma:ap:convolve} to obtain
\begin{align}
  {\mathcal C} &\leq C \Vert u \Vert_{W^{1,\infty}(\Omega)} \Vert u \Vert_{H^3(\Omega)} + C \tau \Vert u \Vert_{L^\infty(\Omega)} \Vert u \Vert_{Y_\tau} \notag\\
  & + C \sum\limits_{m=3}^{\infty} \sum\limits_{j=1}^{m} {m \choose j} \frac{\tau^{m-3}}{(m-3)!^s} \sum\limits_{|\alpha| = m} \sum\limits_{|\beta|=j,\; \beta\leq \alpha} \epsilon^{\alpha_3} \sup\limits_{0<\delta\leq \delta_0} \delta^{m-3} \Vert \partial^\beta u \Vert_{L^p(\Omega_\delta)} \Vert \partial^{\alpha-\beta} Du \Vert_{L^{2p/(p-2)}(\Omega_\delta)},\label{eq:eliminate:theta:5}
\end{align}
where $C$ is a constant depending on the domain, and on $\tau/\tau_* <1$. We rewrite the estimate \eqref{eq:eliminate:theta:5} as
\begin{align}
  {\mathcal C} \leq C \Vert u \Vert_{H^3(\Omega)}^2 + C \Vert u \Vert_{W^{1,\infty}(\Omega)}^2 + C \tau \Vert u \Vert_{L^\infty(\Omega)} \Vert u \Vert_{Y_\tau} + C\left( {\mathcal C}_1 + {\mathcal C}_2 + {\mathcal C}_{\rm low} + {\mathcal C}_{\rm high} + {\mathcal C}_3 + {\mathcal C}_4 + {\mathcal C}_5\right),
\end{align}where for $1\leq j \leq 2$ we denoted
\begin{align}
  {\mathcal C}_1 =& \sum\limits_{m=3}^{\infty} m \sum\limits_{|\alpha|=m} \sum\limits_{|\beta|=1, \beta\leq \alpha} \sup\limits_{0<\delta\leq \delta_0} \Big( \epsilon^{\beta_3} \Vert \partial^\beta u \Vert_{L^\infty(\Omega_\delta)}\Big) \Big( \epsilon^{\alpha_3 - \beta_3} \delta^{m-3} \Vert \partial^{\alpha-\beta} D u \Vert_{L^2(\Omega_\delta)}\Big) \frac{\tau^{m-3}}{(m-3)!^s},\\
  {\mathcal C}_2 =&\sum\limits_{m=4}^{\infty}  {m \choose 2} \sum\limits_{|\alpha|=m} \sum\limits_{|\beta|=2,\beta \leq \alpha} \sup\limits_{0<\delta\leq \delta_0} \Big( \epsilon^{\beta_3} \delta \Vert \partial^\beta u\Vert_{L^\infty(\Omega_{\delta})} \Big)\Big( \epsilon^{\alpha_3 - \beta_3} \delta^{m-4} \Vert \partial^{\alpha-\beta} D u \Vert_{L^2(\Omega_\delta)} \Big)\frac{\tau^{m-3}}{(m-3)!^s},
  \end{align}for $3\leq j \leq m-3$
  \begin{align}
  {\mathcal C}_{\rm low}=&\sum\limits_{m=6}^{\infty} \sum\limits_{j=3}^{[m/2]} {m \choose j} \sum\limits_{|\alpha|=m} \sum\limits_{|\beta|=j,\beta \leq \alpha} \sup\limits_{0<\delta\leq \delta_0} \Big( \epsilon^{\beta_3} \delta^{j-1} \Vert \partial^\beta u\Vert_{L^\infty(\Omega_{\delta})} \Big)\Big( \epsilon^{\alpha_3 - \beta_3} \delta^{m-j-2} \Vert \partial^{\alpha-\beta} D u \Vert_{L^2(\Omega_\delta)} \Big)\frac{\tau^{m-3}}{(m-3)!^s},\\
  {\mathcal C}_{\rm high}=&\sum\limits_{m=7}^{\infty} \sum\limits_{j=[m/2]+1}^{m-3} {m \choose j} \sum\limits_{|\alpha|=m} \sum\limits_{|\beta|=j,\beta \leq \alpha} \sup\limits_{0<\delta\leq \delta_0} \Big( \epsilon^{\beta_3} \delta^{j-3} \Vert \partial^\beta u\Vert_{L^2(\Omega_{\delta})} \Big)\Big( \epsilon^{\alpha_3 - \beta_3} \delta^{m-j} \Vert \partial^{\alpha-\beta} D u \Vert_{L^\infty(\Omega_\delta)} \Big)\frac{\tau^{m-3}}{(m-3)!^s},
  \end{align}and for $m-2\leq j \leq m$
  \begin{align}
  {\mathcal C}_3=&\sum\limits_{m=5}^{\infty}  {m \choose m-2} \sum\limits_{|\alpha|=m} \sum\limits_{|\beta|=m-2,\beta \leq \alpha} \sup\limits_{0<\delta\leq \delta_0} \Big( \epsilon^{\beta_3} \delta^{m-5}\Vert \partial^\beta u\Vert_{L^2(\Omega_{\delta})} \Big)\Big( \epsilon^{\alpha_3 - \beta_3} \delta^2 \Vert \partial^{\alpha-\beta} D u \Vert_{L^\infty(\Omega_\delta)} \Big)\frac{\tau^{m-3}}{(m-3)!^s},\\
  {\mathcal C}_4=&\sum\limits_{m=4}^{\infty}  {m \choose m-1} \sum\limits_{|\alpha|=m} \sum\limits_{|\beta|=m-1,\beta \leq \alpha} \sup\limits_{0<\delta\leq \delta_0} \Big( \epsilon^{\beta_3} \delta^{m-4}\Vert \partial^\beta u\Vert_{L^2(\Omega_{\delta})} \Big)\Big( \epsilon^{\alpha_3 - \beta_3} \delta \Vert \partial^{\alpha-\beta} D u \Vert_{L^\infty(\Omega_\delta)} \Big)\frac{\tau^{m-3}}{(m-3)!^s},\\
  {\mathcal C}_5=&\sum\limits_{m=3}^{\infty} \sum\limits_{|\alpha|=m} \sup\limits_{0<\delta\leq \delta_0} \Big(\epsilon^{\alpha_3} \delta^{m-3} \Vert \partial^\alpha u \Vert_{L^2(\Omega_\delta)}\Big) \Vert D u \Vert_{L^\infty(\Omega_\delta)} \frac{\tau^{m-3}}{(m-3)!^s}.
\end{align}
These seven terms are bounded as in the proof of \cite[Lemma 3.2]{KV2}. Namely, letting $\Omega = \bigcup_{0<\delta\leq \delta_0}\Omega_\delta$, and $j=|\beta|$, for the cases $j=1$ and $j=m$ we have
\begin{align}\label{eq:C1:bound}
  {\mathcal C}_1 + {\mathcal C}_5 \leq C \Vert D u \Vert_{L^\infty(\Omega)} \Vert u \Vert_{H^3(\Omega)} + C \tau \Vert D u \Vert_{L^\infty(\Omega)} \Vert u \Vert_{Y_\tau},
\end{align}for the cases $j=2$ and $j=m-1$ it holds that
\begin{align}
  {\mathcal C}_2 + {\mathcal C}_4 \leq C \tau \Vert D^2 u \Vert_{L^\infty(\Omega)} \Vert u \Vert_{H^3(\Omega)} + C \tau^2 \Vert D^2 u \Vert_{L^\infty(\Omega)} \Vert u \Vert_{Y_\tau},
\end{align}when $j=m-2$ we have
\begin{align}
  {\mathcal C}_3 \leq C \tau^2 \Vert u \Vert_{H^5(\Omega)}^2 + \tau^3 \Vert u \Vert_{H^5(\Omega)} \Vert u \Vert_{Y_\tau},
\end{align}
and when $3\leq j \leq m-2$ we have
\begin{align}
  {\mathcal C}_{\rm low} + {\mathcal C}_{\rm high} &\leq C \left(\tau^{3/2}+(1+K^3) \tau^3\right) \Vert u \Vert_{X_\tau} \Vert u \Vert_{Y_\tau}, \label{eq:Clow+Chigh}
\end{align} for some sufficiently large constant $C$, where $K$ is as defined in \eqref{eq:K:def}. We sketch the proof of the ${\mathcal C}_{\rm low}$ estimate and refer the reader to \cite{KV2} for further details on the other five terms. Modulo multiplying by a smooth cut-off function $\eta$ supported on $\Omega_{\delta-r}$ and which is identically $1$ on $\Omega_\delta$, \eqref{eq:cutoff:estimate} and the three-dimensional Agmon inequality give that
\begin{align}
  \Vert v \Vert_{L^\infty(\Omega_{\delta})} \leq C \Vert v \Vert_{L^2(\Omega_{\delta-r})}^{1/4} \Vert \Delta v \Vert_{L^2(\Omega_{\delta-r})}^{3/4} + \frac{C + C K^{3}}{r^{3/2}} \Vert v \Vert_{L^2(\Omega_{\delta-r})},
\end{align}where $C>0$ is a constant depending on the geometry of $\Omega$ and on $\delta_0$,
and $K(t) = \Vert  u \Vert_{L_t^1(0,t) W_x^{1,\infty}(D)}$ is as in \eqref{eq:K:def}. Letting $r=\delta/j$, for $j\geq 3$, we have
\begin{align}
  & \sup\limits_{0<\delta\leq \delta_0} \Big( \epsilon^{\beta_3} \delta^{j-1} \Vert \partial^\beta u\Vert_{L^\infty(\Omega_{\delta})} \Big) \leq C (1+K^3) \sup\limits_{0<\delta\leq \delta_0} (j/\delta)^{3/2}   \Big( \epsilon^{\beta_3} (\delta-\delta/j)^{j-3} \Vert \partial^\beta u\Vert_{L^2(\Omega_{\delta-\delta/j})} \Big) \delta^2 \notag\\
  & + \sup\limits_{0<\delta\leq \delta_0} \Big( \epsilon^{\beta_3} (\delta-\delta/j)^{j-3} \Vert \partial^\beta u\Vert_{L^2(\Omega_{\delta-\delta/j})} \Big)^{1/4} \sup\limits_{0<\delta\leq \delta_0} \Big( \epsilon^{\beta_3} (\delta-\delta/j)^{j-1} \Vert \partial^\beta u\Vert_{L^2(\Omega_{\delta-\delta/j})} \Big)^{3/4} \delta^{1/2}
\end{align}
In the above inequality we used $(1 + 1/(j-1))^{j-1} \leq e$ for all $j\geq 1$. By the H\"older inequality, and \cite[Lemma 4.2]{KV2}, we obtain from the definition of ${\mathcal C}_{\rm low}$ and the above inequality
\begin{align}
  {\mathcal C}_{\rm low} &\leq C \sum\limits_{m=6}^{\infty} \sum\limits_{j=3}^{[m/2]} {m \choose j} [u]_{j}^{1/4} [u]_{j+2}^{3/4} [u]_{m-j+1} \frac{\tau^{m-3}}{(m-3)!^s}\notag\\
   & \qquad \qquad + C (1+K^3)\sum\limits_{m=6}^{\infty} \sum\limits_{j=3}^{[m/2]} {m \choose j} [u]_{j} j^{3/2} [u]_{m-j+1} \frac{\tau^{m-3}}{(m-3)!^s},\label{eq:Clow:*}
\end{align}where we have denoted
\begin{align}
  [v]_m = \sum\limits_{|\alpha|=m} \epsilon^{\alpha_3} \sup\limits_{0<\delta\leq \delta_0} \delta^{m-3} \Vert \partial^\alpha v \Vert_{L^2(\Omega_\delta)}
\end{align}for all smooth $v$.
In the above estimate \eqref{eq:Clow:*} we used the fact that $[D v]_m \leq C [v]_{m+1}$ and $[\Delta v]_m \leq C [v]_{m+2}$, where $C>0$ may depend on $\epsilon$, which is fixed. The right side of \eqref{eq:Clow:*} is then bounded by
\begin{align*}
  C \left(\tau^{3/2}+(1+K^3) \tau^3\right) \Vert u \Vert_{X_\tau} \Vert u \Vert_{Y_\tau}.
\end{align*}Here we used the definitions \eqref{eq:Xtaudef}--\eqref{eq:Ytaudef}, the discrete Young and H\"older inequalities, and the combinatorial estimate
\begin{align}
  {m \choose j} \frac{(j-3)!^{s/4} (j-1)!^{3s/4} (m-j-2)!^s}{(m-3)!^s (m-j+1)} + {m \choose j} \frac{(j-3)!^s j^{3/2} (m-j-2)!^s}{(m-3)!^s (m-j+1)} \leq C,
\end{align}which holds for all $m\geq6$, $3\leq j \leq [m/2]$, and $s\geq1$, where $C>0$ is a dimensional constant. By reversing the roles of $j$ and $m-j$, similar estimates give the bound on ${\mathcal C}_{high}$, thereby proving \eqref{eq:Clow+Chigh}. This concludes the proof of the Lemma~\ref{lemma:velocitycommutator}.
\end{proof}

\section{The pressure estimate}\label{sec:pressure}\setcounter{equation}{0}
The goal of this section is to prove Lemma~\ref{lemma:pressureterm}. This is achieved in several steps: First we use an $H^2$ regularity estimate on the flattened domain to estimate all tangential derivatives of the pressure; next, we obtain a recursion formula to bootstrap to an estimate with higher number of normal derivatives, which leads to an estimate in terms of the velocity; lastly, we prove a product-type estimate for the Lagrangian Gevrey-class norms defined in Section~\ref{sec:notation} which concludes the proof.

For the rest of this section all functions depend on $y=\theta(x)$, hence we shall further suppress all tildes, and since there is no time evolution for the pressure we also suppress the time dependence.

\subsection*{Semi-norms and a decomposition of the pressure term}
The following semi-norms are useful when treating the pressure term. Namely, define
\begin{align}
  \langle v \rangle_{l,\gamma,\delta} = \delta^{l + |\gamma| - 3} \Vert \partial_1^{\gamma_1} \partial_2^{\gamma_2} v \Vert_{L^2(\Omega_\delta)}\label{eq:<>alpha}
\end{align}for all $\gamma \in {\mathbb N}_{0}^{2}$, $l\in{\mathbb Z}$ with $|\gamma| + l \geq 3$, and all $0<\delta\leq \delta_0$. Also let
\begin{align}
  \langle v\rangle_{l,n} = \sum\limits_{|\alpha| = n,\, \alpha_3 = 0} \sup\limits_{0<\delta\leq \delta_0} \langle v\rangle_{l,\alpha',\delta} = \sum\limits_{|\alpha| = n,\, \alpha_3 = 0} \sup\limits_{0<\delta\leq \delta_0}
  \delta^{l+n-3} \Vert \partial_{1}^{\alpha_1} \partial_2^{\alpha_2} v \Vert_{L^2(\Omega_\delta)}\label{eq:<>n}
\end{align}for all $n\geq 0$ with $n+l \geq 3$. Note that we have the inequality $\langle {D'}^k v\rangle_{l,n} \leq \langle v\rangle_{l-k,n+k}$, where $\langle D'^k v \rangle_{l,n} = \sum\limits_{|\alpha|=n,\alpha_3=0} \sup_{0<\delta\leq \delta_0} \delta^{l+n-3} \Vert D'^k v \Vert_{L^2(\Omega_\delta)}$, and $\Vert D'^k v \Vert_{L^2}$ is defined above.

Next, we shall estimate the pressure term arising in \eqref{eq:ode2}, i.e.,
\begin{align}
  \mathcal P = \sum\limits_{m=3}^{\infty} \left( \sum\limits_{|\alpha|=m} \epsilon^{\alpha_3} \sup\limits_{0<\delta\leq \delta_0} \delta^{m-3} \left\Vert \partial^\alpha \bigl(\partial_j \theta_k\; \partial_k p\bigr) \right\Vert_{L^2(\Omega_\delta)} \right) \frac{\tau^{m-3}}{(m-3)!^s}.
\end{align}Similarly to \eqref{eq:eliminate:theta:1}--\eqref{eq:eliminate:theta:5}, we let $C,\tau_*>0$ be such that $\sum_{|\beta|=j} \Vert  \partial^\beta D \theta \Vert_{L^\infty(\Omega)} \leq C (j-3)!^s / \tau_{*}^{j}$ for all $j\geq 0$. Assuming that $\tau<\tau_*$, it follows from the Leibniz rule and the bound ${m \choose j} (m-j-3)!^s (j-3)!^s (m-3)!^{-s} \leq C$ that the pressure term is bounded as
\begin{align}
  \mathcal P \leq C \Vert D p \Vert_{W^{2,\infty}(\Omega)} + C \sum\limits_{m=3}^{\infty}\left( \sum\limits_{|\alpha|=m}\epsilon^{\alpha_3} \sup\limits_{0<\delta\leq\delta_0} \delta^{m-3} \Vert \partial^\alpha D p \Vert_{L^2(\Omega_\delta)}\right) \frac{\tau^{m-3}}{(m-3)!^s},
\end{align}for some positive constant $C = C(D,\eta)$, where $\eta = \tau/\tau_* < 1$ by assumption. We decompose the upper bound on the pressure term as follows
\begin{align}
 {\mathcal P}& \leq C \Vert D p \Vert_{W^{2,\infty}(\Omega)}+ C \sum\limits_{m=3}^{\infty} \left(\sum\limits_{\alpha_3 = 0}^{m} \epsilon^{\alpha_3} \langle \partial_{3}^{\alpha_3} D p \rangle_{\alpha_3,m-\alpha_3}\right) \frac{\tau^{m-3}}{(m-3)!^s}\notag\\
 & \leq C \Vert p \Vert_{W^{3,\infty}(\Omega)} + C \Bigl( (1+\epsilon){\mathcal P}_0 + {\mathcal P}_1 +  {\mathcal P}_2\Bigr),\label{eq:todo1}
 \end{align} where we have denoted the term with at most one normal derivative by
 \begin{align}
 &{\mathcal P}_0 = \sum\limits_{m=3}^{\infty} \langle D p\rangle_{0,m}  \frac{\tau^{m-3}}{(m-3)!^s}, \label{eq:P0def}
 \end{align}and the terms with at least two normal derivatives (according to $D = \partial_3 + D'$) by
 \begin{align}
 &{\mathcal P}_1 = \sum\limits_{m=3}^{\infty} \sum\limits_{\alpha_3 = 1}^{m} \epsilon^{\alpha_3} \langle \partial_{3}^{\alpha_3 + 1} p\rangle_{\alpha_3,m-\alpha_3}  \frac{\tau^{m-3}}{(m-3)!^s} \label{eq:P1def}
 \end{align} and
 \begin{align}
 &{\mathcal P}_2 = \sum\limits_{m=3}^{\infty} \sum\limits_{\alpha_3 = 2}^{m} \epsilon^{\alpha_3} \langle \partial_{3}^{\alpha_3} p\rangle_{\alpha_3-1,m-\alpha_3+1}  \frac{\tau^{m-3}}{(m-3)!^s} \label{eq:P1def}.
\end{align}
\subsection*{The elliptic Neumann problem for the pressure}
Under the change of variables $\theta\colon x \mapsto y$, the elliptic Neumann problem for the pressure (cf.~\cite{T}) becomes (omitting tildes)
\begin{align}
  & - \Delta p = A_{ij}\ \partial_{ij} p + B_j\ \partial_j p + D_{ijkl}\ \partial_i u_j\  \partial_k u_l,\ \mbox{in}\ \Omega \\
  & - \partial_3 p = C_j\ \partial_j p + \Phi_{ij}\ u_i u_j, \ \mbox{on}\ \partial\Omega,
\end{align}where we denoted
\begin{align}
  & A_{ij} = \frac{1}{\Gamma} \left(
               \begin{array}{ccc}
                 -(\partial_1 \gamma)^2 - (\partial_2 \gamma)^2 & 0 & -\partial_1 \gamma \\
                 0 & -(\partial_1 \gamma)^2 - (\partial_2 \gamma)^2 & -\partial_2 \gamma \\
                 -\partial_1 \gamma & -\partial_2 \gamma & 0 \\
               \end{array}
             \right),\label{eq:Adef}\\
  & B_j = \frac{1}{\Gamma} \left(
                                                                      \begin{array}{c}
                                                                        0 \\
                                                                        0 \\
                                                                        -\partial_{11} \gamma - \partial_{22} \gamma  \\
                                                                      \end{array}
                                                                    \right),\label{eq:Bdef}\\
  & C_j = \frac{1 + \Gamma^{1/2}}{(2 -\Gamma +\Gamma^{1/2})\Gamma^{1/2}} \left(
                \begin{array}{c}
                  \partial_1 \gamma \\
                  \partial_2 \gamma  \\
                  0 \\
                \end{array}
              \right),\label{eq:Cdef}\\
  & D_{ijkl} = \frac{1}{\Gamma} \delta_{ik} \delta_{jl} + \frac{1}{\Gamma}\delta_{k3} \left(
                                                       \begin{array}{ccc}
                                                         \partial_l \gamma& 0 & 0 \\
                                                         0 & \partial_l \gamma & 0 \\
                                                         (\partial_l \gamma)^2 & \partial_1 \gamma \partial_2 \gamma& \partial_l \gamma \\
                                                       \end{array}
                                                     \right),\label{eq:Ddef} \\
  & \Phi_{ij} = \frac{1}{\Gamma^{3/2}} \left(
                                         \begin{array}{ccc}
                                           \partial_{11}\gamma & \partial_{12}\gamma & 0 \\
                                           \partial_{12}\gamma & \partial_{22}\gamma & 0 \\
                                           0 & 0 & 0 \\
                                         \end{array}
                                       \right),\label{eq:Phidef}
  \end{align}
  with
  \begin{align}
  & \Gamma = 1+(\partial_1 \gamma)^2 + (\partial_2 \gamma)^2\label{eq:Gammadef}.
\end{align}The precise form of the above matrices is not essential;  what is important for the following arguments is that $A_{33}=C_3 =0$, and that the coefficients $A_{ij},C_j$ are small. We also denote
\begin{align}
  & f = D_{ijkl}\, \partial_i u_j\, \partial_k u_l,\label{eq:fdef}\\
  & g = \Phi_{ij}\, u_i\, u_j.\label{eq:gdef}
\end{align}

\subsection*{The interior $H^2$-regularity estimate}
Let $p$ be the smooth solution of the elliptic Neumann problem
\begin{align}
  & - \Delta p = A_{ij}\ \partial_{ij} p + B_j\ \partial_j p + f,\ \mbox{in}\ \Omega, \label{eq:P.1}\\
  & - \partial_3 p = C_j\ \partial_j p + g, \ \mbox{on}\ \partial\Omega,\label{eq:P.2}
\end{align} where all coefficients are of Gevrey-class $s$.
We have the following interior $H^2$-regularity estimate.
\begin{lemma}\label{lemma:H^2}
There exists a sufficiently small positive dimensional constant $\overline{\varepsilon}$ such that if $A_{33} = C_3 = 0$, $\Vert A_{ij} \Vert_{L^\infty(\overline{\Omega})} \leq \overline{\varepsilon}$, and $\Vert C_j \Vert_{L^\infty(\overline{\Omega})}\leq \overline{\varepsilon}$, the smooth solution $p$ of \eqref{eq:P.1}--\eqref{eq:P.2} satisfies
\begin{align}
  \Vert D^2 p \Vert_{L^2(\Omega_{\delta+r})} \leq C_0 \left( \Vert f \Vert_{L^2(\Omega_\delta)} + \Vert D g \Vert_{L^2(\Omega_\delta)} + \frac{1+K^2}{r} \Vert D p \Vert_{L^2(\Omega_\delta)}\right),\label{eq:H^2}
\end{align}for all $\delta\in(0,\delta_0)$ and for all $0<r<<1$, where $K$ is as defined in \eqref{eq:K:def}, and $C_0 = C_0(A_{ij},B_j,C_j)$ is a positive constant depending on the domain.
\end{lemma}The proof is standard and thus omitted. It relies on the fact that the elliptic operator acting on $p$ in \eqref{eq:P.1} is a small/lower-order perturbation of the Laplacian, and on the the fact that by \eqref{eq:cutoff:estimate} we have $C_0 {\rm dist}(\Omega_{\delta+r}, \Omega_{\delta}^c) \geq r/(1+K^2)$, for some sufficiently large constant $C_0$.
\subsection*{The estimation of tangential derivatives}
Fix $k\geq 3$, and let $\alpha' = (\alpha_1,\alpha_2,0)\in {\mathbb N}_{0}^{3}$ be such that $|\alpha'| = k$. Consider the system \eqref{eq:P.1}--\eqref{eq:P.2}. The function $\partial^{\alpha'} p$ satisfies the elliptic Neumann problem
\begin{align}
& - \Delta (\partial^{\alpha'} p) = A_{ij}\ \partial_{ij}\partial^{\alpha'}p + B_j\ \partial_j \partial^{\alpha'} p + \partial^{\alpha'} f + [\partial^{\alpha'},A_{ij}\ \partial_{ij}] p + [\partial^{\alpha'}, B_j\ \partial_j] p,\ \mbox{in}\ \Omega,\label{eq:P:tan:1}\\
& - \partial_3 (\partial^{\alpha'} p) = C_j\ \partial_j \partial^{\alpha'} p + \partial^{\alpha'} g + [\partial^{\alpha'},C_j\ \partial_j] p,\ \mbox{on}\ \partial\Omega\label{eq:P:tan:2}.
\end{align}We apply the $H^2$-estimate \eqref{eq:H^2} to the solution of \eqref{eq:P:tan:1}--\eqref{eq:P:tan:2}, and bound the commutators using the Leibniz rule as
\begin{align}
  \Vert[\partial^{\alpha'},A_{ij}\ \partial_{ij}] p\Vert_{L^2(\Omega_\delta)} \leq \sum\limits_{0< \beta' \leq \alpha'} {\alpha' \choose \beta'} \Vert \partial^{\beta'} A_{ij} \Vert_{L^\infty(\Omega_\delta)} \Vert \partial^{\alpha'-\beta'} \partial_{ij} p \Vert_{L^2(\Omega_\delta)}.
\end{align}The terms involving $[\partial^{\alpha'},B_j \partial_j]$ and $[\partial^{\alpha'},C_j \partial_j]$ are estimated similarly. Letting $r=\delta/k \leq \delta/3$, we obtain
\begin{align}
  \Vert \partial^{\alpha'} D^2 p \Vert_{L^2(\Omega_{\delta + \delta/k})} &\leq  C_0 \left( \Vert \partial^{\alpha'}f \Vert_{L^2(\Omega_\delta)} + \Vert \partial^{\alpha'} D g \Vert_{L^2(\Omega_\delta)} + (1+K^2) \frac{k}{\delta} \Vert \partial^{\alpha'} D p \Vert_{\Omega_\delta}\right)\notag \\
  & + C_0 \sum\limits_{j=1}^{k} {k \choose j} \sum\limits_{|\beta'| =j,\ \beta'\leq \alpha'} \max\{\Vert \partial^{\beta'} A_{ij} \Vert_{L^\infty(\bar{\Omega}_\delta)}, \Vert \partial^{\beta'} B_j \Vert_{L^\infty(\bar{\Omega}_\delta)}, \Vert \partial^{\beta'} C_j\Vert_{L^\infty(\bar{\Omega}_\delta)} \}\notag\\
   & \qquad \qquad \qquad \qquad \qquad \qquad \times \left( \Vert\partial^{\alpha'-\beta'} D D'  p \Vert_{L^2(\Omega_\delta)} + \Vert \partial^{\alpha'-\beta'} D p \Vert_{L^2(\Omega_\delta)} \right), \label{eq:rough:H2}
\end{align}where we used $A_{33}=0$. Denote
\begin{align}
\psi_{\beta',\delta}=  \delta^{|\beta'|-2} \max\{\Vert \partial^{\beta'} A_{ij} \Vert_{L^\infty(\bar{\Omega}_\delta)}, \Vert \partial^{\beta'} B_j \Vert_{L^\infty(\bar{\Omega}_\delta)}, \Vert \partial^{\beta'} C_j\Vert_{L^\infty(\bar{\Omega}_\delta)} \},\label{eq:psi:beta:def}
 \end{align}for all $|\beta'|\geq 2$, and
 \begin{align}
\psi_{\beta',\delta}=  \max\{\Vert \partial^{\beta'} A_{ij} \Vert_{L^\infty(\bar{\Omega}_\delta)}, \Vert \partial^{\beta'} B_j \Vert_{L^\infty(\bar{\Omega}_\delta)}, \Vert \partial^{\beta'} C_j\Vert_{L^\infty(\bar{\Omega}_\delta)} \},\label{eq:psi:beta:low}
 \end{align} if $|\beta'|=1$. Also let
\begin{align}
  \psi_j = \sum\limits_{|\beta'|=j} \sup\limits_{0<\delta\leq \delta_0} \psi_{\beta',\delta},\label{eq:psi:def}
\end{align}Note that since $\gamma$, and hence $A_{ij},B_j, C_j$, is of Gevrey-class $s$, there exist $C,\tau_*>0$ such that \begin{align}
    \psi_j \leq C \frac{(j-2)!^s}{\tau_{*}^{j}},\label{eq:psi:bound}
  \end{align}for all $j\geq1$, where recall $(-1)!=1$. Since the boundary is fixed, $C$ and the Gevrey-class $s$ radius $\tau_*$ of the boundary are not functions of time. Multiplying estimate \eqref{eq:rough:H2} by $(\delta + \delta/k)^{k-2}$, and using $(1+ 1/k)^{k-2} \leq e$ for $k\geq 3$, we obtain
\begin{align}
  &\langle D^2 p\rangle_{1,\alpha',\delta + \delta/k} \leq C \langle f\rangle_{1,\alpha',\delta} + C  \langle D g\rangle_{1,\alpha',\delta} + C \psi_{\alpha',\delta} \left( \Vert DD' p \Vert_{L^2(\Omega_\delta)}+\Vert Dp \Vert_{L^2(\Omega_\delta)}\right)\notag\\
  &\qquad \qquad + C (1+K^2) k \sum\limits_{|\beta'|=k-1,\beta'\leq \alpha'} \psi_{\beta',\delta} \Big( \delta \Vert D D'^2 p\Vert_{L^2(\Omega_\delta)} + \delta \Vert D D'p\Vert_{L^2(\Omega_\delta)}\Big)\notag\\
  &\qquad  \qquad + C k (k-1)\sum\limits_{|\beta'|=k-2,\beta'\leq \alpha'} \psi_{\beta',\delta} \delta \Vert D D'^2 p\Vert_{L^2(\Omega_\delta)}\notag\\
  &\qquad \qquad + C \chi_{k\geq 4} \sum\limits_{j=2}^{k-2} {k \choose j}  \sum\limits_{|\beta'|=j, \beta'\leq \alpha'}  \psi_{\beta',\delta} \langle D D' p\rangle_{3,\alpha'-\beta',\delta} + C \chi_{k\geq 5} \sum\limits_{j=2}^{k-3} {k \choose j} \sum\limits_{|\beta'|=j, \beta'\leq \alpha'} \psi_{\beta',\delta} \langle D p\rangle_{3,\alpha'-\beta',\delta} \notag\\
  &\qquad \qquad +C k \langle D p\rangle_{0,\alpha',\delta} + C k \sum\limits_{|\beta'|=1,\beta'\leq \alpha'} \psi_{\beta',\delta} \Big( \langle D D' p\rangle_{2,\alpha'-\beta',\delta} + \langle Dp\rangle_{2,\alpha'-\beta',\delta}\Big) .\label{eq:deltarole}
\end{align}By taking the supremum over $0<\delta\leq \delta_0 \leq 1$ of the above estimate and summing over all $|\alpha'|=k\geq 3$, cf.~\cite[Lemma 4.2]{KV2}, we obtain
\begin{align}
  \langle \partial_3 Dp\rangle_{1,k} &+ \langle Dp \rangle_{0,k+1}\leq C \Bigl( \langle f \rangle_{0,k} +\langle D g\rangle_{0,k}\Bigr) \notag\\
  &\qquad \qquad  + C (1+K^2) \Bigl( \psi_k \Vert Dp \Vert_{L^2(\Omega)}  + (\psi_k + k \psi_{k-1}) \Vert D^2 p \Vert_{L^2(\Omega)} + (k \psi_{k-1} + k^2 \psi_{k-2} )\Vert D^3 p \Vert_{L^2(\Omega)} \Bigr)\notag\\
  &\qquad \qquad +C \sum\limits_{j=1}^{k-2} {k \choose j} \psi_j \langle Dp\rangle_{0,k-j+1} +C \chi_{k\geq 4} \sum\limits_{j=1}^{k-3} {k \choose j} \psi_j \langle Dp\rangle_{0,k-j},\label{eq:tangentialH^2}
\end{align} where as usual we write $\Omega = \bigcup_{0<\delta\leq \delta_0} \Omega_\delta$. Estimate \eqref{eq:tangentialH^2} is used to bound the term ${\mathcal P}_0$ in the decomposition \eqref{eq:todo1} of ${\mathcal P}$. Furthermore, using the bound \eqref{eq:psi:bound} on $\psi_j$, estimate \eqref{eq:todo1} implies
\begin{align}
\langle \partial_3 Dp\rangle_{1,k} + \langle Dp \rangle_{0,k+1} &\leq C \Bigl( \langle f\rangle_{0,k} +\langle D g\rangle_{0,k}\Bigr) \notag\\
  & + C  (1+\tau_* + \tau_{*}^{2})(1+K^2) \Bigl( \Vert Dp \Vert_{L^2(\Omega)}  + \Vert D^2 p \Vert_{L^2(\Omega)} + \Vert D^3 p \Vert_{L^2(\Omega)} \Bigr) \frac{(k-2)!^s}{\tau_{*}^k}\notag\\
  &+C  \sum\limits_{j=1}^{k-2} {k \choose j} \frac{(j-2)!^s}{\tau_{*}^j} \langle Dp\rangle_{0,k-j+1} +C \chi_{k\geq 4}  \sum\limits_{j=1}^{k-3} {k \choose j} \frac{(j-2)!^s}{\tau_{*}^j} \langle Dp\rangle_{0,k-j},\label{eq:tangentialH^2:new}
\end{align}for all $k\geq 3$, where $C$ depends on $C_0$ and $\delta_0\leq 1$, while $\tau_*$ is fixed, depending only on $\gamma$.
\subsection*{The transfer of normal to tangential derivatives}
We use the special structure of the coefficients $A_{ij}$ and $B_j$ to rewrite \eqref{eq:P.1} as
\begin{align}
  -\partial_{33} p &= (a_1\ \partial_1 + a_2\ \partial_2) \partial_3 p + b\ \partial_3 p + c\ (\partial_{11} + \partial_{22}) p + f \notag\\
   &= (a\cdot\nabla') \partial_3 p + b\ \partial_3 p + c\ \Delta' p + f \label{eq:P1new},
\end{align}where, as above, (cf.~\eqref{eq:Adef},\eqref{eq:Bdef}, and \eqref{eq:Gammadef})
\begin{align}
a_i = - 2\frac{\partial_i \gamma}{\Gamma},\ b = - \frac{\partial_{11} \gamma + \partial_{22} \gamma}{\Gamma},\ c =\frac{1}{\Gamma}.
\end{align}Since $\gamma$, and hence $a,b$, and $c$, is a function of $(y_1,y_2)$ only, we obtain from \eqref{eq:P1new} that for $k \geq 2$ we have
\begin{align}
-  \partial_{3}^{k} p = (a\cdot\nabla') \partial_{3}^{k-1} p + b\ \partial_{3}^{k-1} p + c\ \Delta' \partial_{3}^{k-2} p + \partial_{3}^{k-2} f.\label{eq:D3^k}
\end{align}Note that in the case of the half-space (cf.~\cite{KV2}), identity \eqref{eq:P1new} simplifies to $-\partial_{33} p = \Delta' p + f$, which allows one to obtain an explicit formula for $\partial_{33}^k p$ in terms of $f$ and $(-\Delta')^k p$. The combinatorial structure of this transfer of normal to tangential derivatives is encoded in the coefficients $M_\alpha$ of \cite{KV2}. In the case of the present paper, it is highly inconvenient use the recursion formula \eqref{eq:D3^k} to explicitly calculate $\partial_{3}^k p$ in terms of $f$ and tangential derivatives of $p$. Instead we use the fact that we may choose $\epsilon <<1$ and recursively bootstrap to estimates on higher number of normal derivatives acting on $p$. By applying $\partial^{\alpha'}$, where $|\alpha'|=n$, to \eqref{eq:D3^k}, using the Leibniz rule, the H\"older inequality, we obtain
\begin{align}
  & \Vert \partial_{3}^{k} \partial^{\alpha'} p \Vert_{L^2(\Omega_\delta)} \leq \Vert \partial_{3}^{k-2} \partial^{\alpha'} f \Vert_{L^2(\Omega_\delta)}\notag\\
  & + C \sum\limits_{j=0}^{n} {n \choose j} \sum\limits_{|\beta'|=j,\beta'\leq \alpha'} \psi_{\beta',\delta} \left( \Vert \partial^{\alpha'-\beta'}D' \partial_{3}^{k-1} p\Vert_{L^2(\Omega_\delta)}+ \Vert \partial^{\alpha'-\beta'}\partial_{3}^{k-1} p \Vert_{L^2(\Omega_\delta)} + \Vert\partial^{\alpha'-\beta'}\Delta' \partial_{3}^{k-2}p \Vert_{L^2(\Omega_\delta)}\right),\label{eq:D3:ancient}
\end{align}where we have denoted $\psi_{\beta',\delta}$ and $\psi_j$ similarly to \eqref{eq:psi:beta:def}--\eqref{eq:psi:def} (replace $A_{ij},B_j,C_j$ by $a,b,c$). Since $a,b$, and $c$ are of Gevrey-class $s$ (they only depend on $\gamma$), as in  \eqref{eq:psi:bound}, there exist $C,\tau_* >0$ with $\psi_j \leq C (j-2)!^s/\tau_{*}^{j}$. Multiplying the bound \eqref{eq:D3:ancient} by $\delta^{n+k-4}$, it follows that
\begin{align}
  &\langle \partial_{3}^{k} p\rangle_{k-1,\alpha',\delta} \leq \langle \partial_{3}^{k-2} f\rangle_{k-1,\alpha',\delta} \notag\\
  &\ \ \ \ + \sum\limits_{j=0}^{n} {n \choose j} \sum\limits_{|\beta'|=j,\beta'\leq \alpha'}\psi_{\beta',\delta} \Big(  \langle \partial_{3}^{k-1} D' p\rangle_{k+1,\alpha'-\beta',\delta} + \langle\partial_{3}^{k-2} \Delta' p\rangle_{k+1,\alpha'-\beta',\delta} + \langle\partial_{3}^{k-1} p\rangle_{k+1,\alpha'-\beta',\delta} \Big), \label{eq:highD3}
\end{align}for all $n+k \geq 4$, $0\leq \delta\leq \delta_0$, and  $\alpha\in {\mathbb N}_{0}^3$ with $\alpha_3 = 0$, and $|\alpha'|=n$. Estimate \eqref{eq:highD3} above will be used to bound the terms with high number of normal derivatives in the pressure estimate, namely ${\mathcal P}_1$, and ${\mathcal P}_2$.

\subsection*{Bounds for ${\mathcal P}_0, {\mathcal P}_1$, and ${\mathcal P}_2$}
For the term ${\mathcal P}_0$ with a low number of tangential derivatives we have the bound
\begin{align}
  {\mathcal P}_0 &\leq \frac{C \eta}{1-\eta} {\mathcal P}_0 + C_1 (1+K^2) \Vert p \Vert_{H^4(\Omega)}+ C_1 \sum\limits_{m=4}^{\infty} \Bigl(\langle f\rangle_{0,m-1} + \langle Dg\rangle_{0,m-1}\Bigr) \frac{\tau^{m-3}}{(m-3)!^s},\label{eq:P0:bound}
\end{align}where $\eta = \tau/\tau_* < 1$, $\tau_*$ is the Gevrey-class radius of the boundary, $C=C(\gamma)$ and $C_1=C_1(\gamma,\eta,\epsilon)$ are sufficiently large constant positive constants. As usual, $\Omega = \bigcup_{0<\delta\leq \delta_0} \Omega_\delta$. Note that the condition $\eta<1$ is natural, as the flow may not have arbitrarily large radius of Gevrey-class regularity close to the boundary. Under the assumption $\eta < 1$, we also have the bounds
\begin{align}
  {\mathcal P}_1 \leq \frac{\epsilon\, C}{1-\eta} \left({\mathcal P}_0 +  {\mathcal P}_1 \right) + C_1 \Vert p \Vert_{H^3(\Omega)} + C_1 \sum\limits_{m=3}^{\infty} \sum\limits_{\alpha_3 = 1}^{m} \epsilon^{\alpha_3} \langle \partial_{3}^{\alpha_3-1} f\rangle_{\alpha_3,m-\alpha_3} \frac{\tau^{m-3}}{(m-3)!^s},\label{eq:P1:bound}
\end{align}and
\begin{align}
  {\mathcal P}_2 \leq &+ \frac{\epsilon\, C}{1-\eta} \left(\epsilon\, {\mathcal P}_0 + {\mathcal P}_1 +  {\mathcal P}_2 \right)+ C_1 \Vert p \Vert_{H^3(\Omega)}  + C_1 \sum\limits_{m=3}^{\infty} \sum\limits_{\alpha_3 = 2}^{m} \epsilon^{\alpha_3} \langle \partial_{3}^{\alpha_3 - 2} f\rangle_{\alpha_3-1,m-\alpha_3+1} \frac{\tau^{m-3}}{(m-3)!^s} ,\label{eq:P2:bound}
\end{align}where $C=C(\gamma)$ is a fixed sufficiently large constant, while $C_1 = C_1(\gamma,\eta,\epsilon)$ has additional dependence on the Gevrey-class norm and radius of $\gamma$ cf.~\eqref{eq:psi:bound} and the parameter $\epsilon$. First we prove the bound for ${\mathcal P}_0$.
\begin{proof}[Proof of \eqref{eq:P0:bound}]
Letting $k=m-1$ in \eqref{eq:tangentialH^2:new}, and recalling the definition \eqref{eq:P0def} of ${\mathcal P}_0$ , we obtain
\begin{align*}
  {\mathcal P}_0 \leq \langle Dp\rangle_{0,3} & + C \sum\limits_{m=4}^{\infty} \Bigl(\langle f\rangle_{0,m-1} + \langle Dg\rangle_{0,m-1} \Bigr) \frac{\tau^{m-3}}{(m-3)!^s} \notag\\
  & + C (1+K^2)(1+\tau_{*}^{2}) \Vert p \Vert_{H^3(\Omega)} \sum\limits_{m=4}^{\infty} \frac{(m-3)!^s}{\tau_{*}^{m-3}} \frac{\tau^{m-3}}{(m-3)!^s}\notag\\
  & + C \sum\limits_{m=4}^{\infty} \sum\limits_{j=1}^{m-3} {m-1 \choose j} \frac{(j-2)!^s}{\tau_{*}^j} \langle Dp\rangle_{0,m-j} \frac{\tau^{m-3}}{(m-3)!^s}\notag\\
  & + C \sum\limits_{m=5}^{\infty} \sum\limits_{j=1}^{m-4} {m-1 \choose j} \frac{(j-2)!^s}{\tau_{*}^j} \langle Dp\rangle_{0,m-j-1} \frac{\tau^{m-3}}{(m-3)!^s}.
\end{align*}Using the fact that for all $s\geq 1$, $m\geq 4$, and $1\leq j \leq m-3$ we have
\begin{align}
  {m-1 \choose j} \frac{(j-2)!^s (m-j-3)!^s}{(m-3)!^s} \leq C,
\end{align}and recalling that we have $\eta = \tau/\tau_* < 1$, we estimate the discrete convolution and obtain
\begin{align}
  {\mathcal P}_0 &\leq   \frac{\eta\, C}{1-\eta} {\mathcal P}_0 + C_1 (1+K^2) \Vert p \Vert_{H^4(\Omega)} + C_1 \sum\limits_{m=4}^{\infty} \Bigl(\langle f\rangle_{0,m-1} + \langle Dg\rangle_{0,m-1} \Bigr) \frac{\tau^{m-3}}{(m-3)!^s},
\end{align}where $C$ is a dimensional constant and $C_1 = C_1(\gamma,\tau_0,\eta)$, concluding the proof.\end{proof}
The estimates for ${\mathcal P}_1$ and ${\mathcal P}_2$ are symmetric, and so to avoid redundancy we only give the proof of \eqref{eq:P1:bound}.
\begin{proof}[Proof of \eqref{eq:P1:bound}]
Let $k=\alpha_3 + 1$ and $n = m-\alpha_3$ (so that $n+k \geq 4$) in \eqref{eq:highD3}, to obtain that for all $|\alpha'|=m-\alpha_3$ we have
\begin{align*}
\langle \partial_{3}^{\alpha_3+1} p\rangle_{\alpha_3,\alpha',\delta} & \leq \langle \partial_{3}^{\alpha_3 -1} f\rangle_{\alpha_3,\alpha',\delta} + \sum\limits_{j=0}^{m-\alpha_3} {m-\alpha_3 \choose j} \sum\limits_{|\beta'|=j,\beta'\leq \alpha'}\psi_{\beta',\delta}\times \notag\\
&\qquad \qquad \qquad \times \Big( \langle \partial_{3}^{\alpha_3} D' p\rangle_{\alpha_3+2,\alpha'-\beta',\delta} + \langle\partial_{3}^{\alpha_3-1} D'^2 p\rangle_{\alpha_3+2,\alpha'-\beta',\delta} + \langle\partial_{3}^{\alpha_3} p\rangle_{\alpha_3+2,\alpha'-\beta',\delta} \Big).
\end{align*}Taking the supremum over $0<\delta\leq \delta_0 <1$, and summing over all $|\alpha'|=m-\alpha_3$, the above estimate implies
\begin{align*}
  {\mathcal P}_1  \leq C \sum\limits_{m=3}^{\infty} &\sum\limits_{\alpha_3 =1}^{m} \epsilon^{\alpha_3} \langle \partial_{3}^{\alpha_3-1} f\rangle_{\alpha_3,m-\alpha_3} \frac{\tau^{m-3}}{(m-3)!^s}+ C \sum\limits_{m=3}^{\infty} \sum\limits_{\alpha_3=1}^{m}\sum\limits_{j=0}^{m-\alpha_3} {m-\alpha_3 \choose j} \psi_j \epsilon^{\alpha_3} \times \notag\\
  &\qquad \qquad \times  \Big(\langle  \partial_{3}^{\alpha_3}p\rangle_{\alpha_3+1,m-j-\alpha_3+1} + \langle \partial_{3}^{\alpha_3-1}p\rangle_{\alpha_3,m-j-\alpha_3+2} + \langle \partial_{3}^{\alpha_3} p\rangle_{\alpha_3+2,m-j-\alpha_3}\Big) \frac{\tau^{m-3}}{(m-3)!^s}.
\end{align*} Using the bound \eqref{eq:psi:bound} on $\psi_j$ and the combinatorial estimate
\begin{align}
  {m-\alpha_3 \choose j} \frac{(j-2)!^s (m-j-3)!^s}{(m-3)!^s} \leq C,
\end{align}which holds for all $m\geq 3$, $1\leq \alpha_{3} \leq m$, and $0\leq j \leq m-\alpha_3$, we obtain
\begin{align}
  &{\mathcal P}_1 \leq C \sum\limits_{m=3}^{\infty} \sum\limits_{\alpha_3 =1}^{m} \epsilon^{\alpha_3} \langle\partial_{3}^{\alpha_3-1} f\rangle_{\alpha_3,m-\alpha_3} \frac{\tau^{m-3}}{(m-3)!^s}\notag\\
  & \qquad + \epsilon C \sum\limits_{m=3}^{\infty} \sum\limits_{\alpha_3 = 1}^{m}\sum\limits_{j=0}^{m-\alpha_3} \eta^j \left( \epsilon^{\alpha_3-1} \langle \partial_{3}^{\alpha_3} p\rangle_{\alpha_3+1,m-j-\alpha_3+1} \frac{\tau^{m-j-3}}{(m-j-3)!^s}\right)\notag\\
  & \qquad + \tau \epsilon C\sum\limits_{m=3}^{\infty} \sum\limits_{\alpha_3 = 1}^{m}\sum\limits_{j=0}^{m-\alpha_3} \eta^{j+1} \left( \epsilon^{\alpha_3-1} \langle \partial_{3}^{\alpha_3} p\rangle_{\alpha_3+2,m-j-\alpha_3} \frac{\tau^{m-j-4}}{(m-j-4)!^s}\right)\notag\\
  & \qquad + \epsilon^2 C \sum\limits_{m=3}^{\infty} \sum\limits_{\alpha_3 = 1}^{m}\sum\limits_{j=0}^{m-\alpha_3} \eta^j \left( \epsilon^{\alpha_3-2} \langle \partial_{3}^{\alpha_3-1} p\rangle_{\alpha_3,m-j-\alpha_3+2} \frac{\tau^{m-j-3}}{(m-j-3)!^s}\right).
\end{align}Here, as before we denoted $\eta = \tau/\tau_*<1$. It is convenient to reverse the summation order in the above estimate and write
\begin{align}
  &{\mathcal P}_1 \leq C \sum\limits_{m=3}^{\infty} \sum\limits_{\alpha_3 =1}^{m} \epsilon^{\alpha_3} \langle \partial_{3}^{\alpha_3-1} f\rangle_{\alpha_3,m-\alpha_3} \frac{\tau^{m-3}}{(m-3)!^s}\notag\\
  & \qquad + \epsilon C \sum\limits_{m=3}^{\infty} \sum\limits_{j = 0}^{m-1} \eta^j \sum\limits_{\alpha_3=0}^{m-j-1} \left( \epsilon^{\alpha_3} \langle \partial_{3}^{\alpha_3+1} p\rangle_{\alpha_3+2,m-j-\alpha_3} \frac{\tau^{m-j-3}}{(m-j-3)!^s}\right)\notag\\
  & \qquad + \tau \epsilon C \sum\limits_{m=3}^{\infty} \sum\limits_{j= 1}^{m} \eta^j \sum\limits_{\alpha_3+3,\alpha_3=0}^{m-j} \left( \epsilon^{\alpha_3} \langle \partial_{3}^{\alpha_3+1} p\rangle_{m-j-\alpha_3} \frac{\tau^{m-j-3}}{(m-j-3)!^s}\right)\notag\\
  & \qquad + \epsilon^2 C \sum\limits_{m=3}^{\infty} \sum\limits_{j=0}^{m-2} \eta^j \sum\limits_{\alpha_3=0}^{m-j-2} \left( \epsilon^{\alpha_3} \langle \partial_{3}^{\alpha_3+1} p\rangle_{\alpha_3+2,m-j-\alpha_3} \frac{\tau^{m-j-3}}{(m-j-3)!^s}\right)\notag\\
  & \qquad + \epsilon C \sum\limits_{m=3}^{\infty} \sum\limits_{j=0}^{m-1} \eta^j \langle p\rangle_{1,m-j+1} \frac{\tau^{m-j-3}}{(m-j-3)!^s}\notag\\
  & \leq C \sum\limits_{m=3}^{\infty} \sum\limits_{\alpha_3 =1}^{m} \epsilon^{\alpha_3} \langle\partial_{3}^{\alpha_3-1} f\rangle_{\alpha_3,m-\alpha_3} \frac{\tau^{m-3}}{(m-3)!^s} + T_1 + T_2 + T_3 + T_4.\label{eq:*1}
\end{align}The terms $T_1,T_2,T_3$, and $T_4$ are bounded by estimating the discrete convolution $\sum_m \sum_j x_j y_{m-j}$, and using the fact that since $\eta<1$ we have $\sum_{j\geq 0} \eta^j = 1/(1-\eta)$. We have the following estimate
\begin{align}
  T_1 \leq \frac{\epsilon\, C}{1-\eta} {\mathcal P}_0 +  \frac{\epsilon\, C}{1-\eta} {\mathcal P}_1+ C_1 \Vert p\Vert_{H^3(\Omega)},
\end{align}where $C=C(\gamma)$ is a positive constant, and $C_1$ has additional dependance on $\eta$, and $\epsilon$. Similarly we obtain
\begin{align}
  T_2 & \leq  \frac{\epsilon\, C}{1-\eta} {\mathcal P}_0 +  \frac{\epsilon\, C}{1-\eta} {\mathcal P}_1 + C_1 \Vert p\Vert_{H^3(\Omega)},\\
  T_3 & \leq \frac{\epsilon^2 C}{1-\eta} {\mathcal P}_0 +  \frac{\epsilon^2 C}{1-\eta} {\mathcal P}_1 + C_1 \Vert p\Vert_{H^3(\Omega)},
\end{align} and
\begin{align}
  T_4 \leq  \frac{\epsilon\, C}{1-\eta} {\mathcal P}_0 + C_1  \Vert p\Vert_{H^3(\Omega)},\label{eq:*2}
\end{align}with $C=C(\gamma)>0$, and $C_1 = C_1(\gamma,\epsilon,\eta) > 0$. The proof is concluded by combining \eqref{eq:*1}--\eqref{eq:*2}.\end{proof}

\subsection*{Gevrey-class estimates for the pressure}
\begin{lemma}\label{lemma:mathcal:P}
  There exists a sufficiently small constant $\epsilon > 0$ depending only on $\gamma$, such that if $\tau \leq \epsilon \tau_*$, then we have
  \begin{align}
  {\mathcal P} & \leq C \sum\limits_{m=3}^{\infty} \sum\limits_{\alpha_3=1}^{m} \epsilon^{\alpha_3} \langle\partial_{3}^{\alpha_3-1} f\rangle_{\alpha_3,m-\alpha_3} \frac{\tau^{m-3}}{(m-3)!^s} + C \sum\limits_{m=3}^{\infty} \sum\limits_{\alpha_3=2}^{m} \epsilon^{\alpha_3} \langle \partial_{3}^{\alpha_3-2} f\rangle_{\alpha_3-1,m-\alpha_3+1} \frac{\tau^{m-3}}{(m-3)!^s} \notag \\
  & \qquad + C \sum\limits_{m=3}^{\infty} \langle D g\rangle_{1,m-1} \frac{\tau^{m-3}}{(m-3)!^s} + C (1+K^2) \Vert p \Vert_{H^4(\Omega)} + C \Vert p \Vert_{W^{3,\infty}(\Omega)},\label{eq:mathcal:P}
\end{align}where $C = C(\gamma)$ is a fixed positive constant.
\end{lemma}
\begin{proof}[Proof of Lemma~\ref{lemma:mathcal:P}]
By combining estimates \eqref{eq:P0:bound}--\eqref{eq:P2:bound} we obtain that for $\eta<1$
\begin{align}
  {\mathcal P}_{0} + {\mathcal P}_1 &+ {\mathcal P}_2 \leq  \frac{(\epsilon + \eta)C_*}{1-\eta} {\mathcal P}_0 + \frac{\epsilon\, C_*}{1-\eta} {\mathcal P}_1 + \frac{\epsilon\, C_*}{1-\eta} {\mathcal P}_2 +  C_1 (1+K^2) \Vert p \Vert_{H^4(\Omega)}  \notag\\
  &+ C_1 \sum\limits_{m=4}^{\infty} \Bigl( \langle f\rangle_{0,m-1} + \langle Dg\rangle_{0,m-1}\Bigr) \frac{\tau^{m-3}}{(m-3)!^s} + C_1 \sum\limits_{m=3}^{\infty} \sum\limits_{\alpha_3=2}^{m} \epsilon^{\alpha_3} \langle D \partial_{3}^{\alpha_3-2} f\rangle_{\alpha_3,m-\alpha_3} \frac{\tau^{m-3}}{(m-3)!^s},\label{eq:P:estimate}
\end{align}for a sufficiently large, fixed constant $C_*=C_*(\gamma)>0$, and $C_1 = C_1(\gamma,\epsilon,\eta)>0$. Define $\epsilon = \epsilon(\gamma)$ by
\begin{align}
  \epsilon = \frac{1}{ 1 + 4C_*}.\label{eq:epsilon:choice}
\end{align}It is clear that $\epsilon$ may be fixed for all time, as it only depends on the boundary of the domain. Whenever $\tau\leq \epsilon \tau_*$, we have $\eta = \tau/\tau_* \leq \epsilon$, and therefore $(\epsilon + \eta)/(1-\eta) \leq 2\epsilon/(1-\epsilon) \leq 1/ (2 C_*)$, by the choice of $\epsilon$ \eqref{eq:epsilon:choice}. Thus the terms involving ${\mathcal P}_0$, ${\mathcal P}_1$, and ${\mathcal P}_2$ on the right side of \eqref{eq:P:estimate}, may be absorbed on the left side of \eqref{eq:P:estimate} and the proof of the lemma is completed.\end{proof}
\begin{remark}
  The condition $\tau < \epsilon \tau_*$ is not restrictive; it is a manifestation of the fact that the velocity field cannot have arbitrarily large Gevrey-class radius close to the boundary, it must be bounded from above by the Gevrey-class radius of the boundary.
\end{remark}
In the following lemma we use the definitions of $f$ and $g$ (cf.~\eqref{eq:fdef}, \eqref{eq:gdef}) to bound the right side of \eqref{eq:mathcal:P} in terms of the velocity.
\begin{lemma}\label{lemma:P:u}
  For $\epsilon = \epsilon(\gamma)>0$ as in Lemma~\ref{lemma:mathcal:P}, if $\tau < \epsilon \tau_*$, then we have
 \begin{align}
{\mathcal P} &\leq C \sum\limits_{m=3}^{\infty} \langle D( u u)\rangle_{1,m-1} \frac{\tau^{m-3}}{(m-3)!^s} + C \sum\limits_{m=3}^{\infty} \sum\limits_{\alpha_3=0}^{m-1}\epsilon^{\alpha_3+1} \langle \partial_3^{\alpha_3}(Du Du ) \rangle_{\alpha_3+3,m-1-\alpha_3} \frac{\tau^{m-3}}{(m-3)!^s}\notag\\
 & \qquad + C \Vert Du Du \Vert_{H^1(\Omega)} + C \Vert u u \Vert_{H^2(\Omega)} + C (1+K^2) \Vert p \Vert_{H^4(\Omega)} + C \Vert p \Vert_{W^{3,\infty}(\Omega)}, \label{eq:P:u}
 \end{align}where $C = C(\gamma)>0$ is a sufficiently large constant.
\end{lemma}
\begin{proof}[Proof of Lemma~\ref{lemma:P:u}]
Denote the right side of \eqref{eq:mathcal:P} by $T_{f}^1 + T_{f}^2 + T_g + C (1+K^2) \Vert p \Vert_{H^4(\Omega)} + C \Vert p \Vert_{W^{2,\infty}(\Omega)}$. First we estimate the term
\begin{align}
  T_g = C \sum\limits_{m=3}^{\infty} \langle Dg\rangle_{1,m-1} \frac{\tau^{m-3}}{(m-3)!^s}.
\end{align}Recall that $g = u_i u_j \Phi_{ij}$ (cf.~\eqref{eq:gdef}). As in the proof of \eqref{eq:P0:bound} and \eqref{eq:P1:bound}, we denote $\psi_{\beta,\delta} = \delta^{|\beta|} \Vert \partial^{\beta} \Phi_{ij} \Vert_{L^\infty(\Omega_\delta)}$, and $\psi_{j} = \sum_{|\beta|=j} \sup_{0<\delta\leq \delta_0} \psi_{\beta,\delta}$. Since $\Phi_{ij}$ is of Gevrey-class $s$ (cf.~\eqref{eq:Phidef}) there exist $C,\tau_*$ such that $\psi_j \leq C (j-3)!^s/\tau_{0}^{*}$, for all $j\geq 0$ (recall that we write $n!=1$ if $n\leq 0$). By the Leibniz rule, we have
\begin{align}
  T_g \leq C \sum\limits_{m=3}^{\infty} \sum\limits_{j=0}^{m} \sum\limits_{|\alpha|=m,\alpha_3 \leq 1} \sum\limits_{|\beta|=j,\beta\leq\alpha} {\alpha \choose \beta} \sup\limits_{0<\delta\leq\delta_0} \delta^{m-3} \Vert \partial^\beta \Phi \Vert_{L^\infty(\Omega_\delta)} \Vert \partial^{\alpha-\beta} (u u) \Vert_{L^2(\Omega_\delta)} \frac{\tau^{m-3}}{(m-3)!^s}.
\end{align}We split this sum into four pieces according to $j=m$, $j=m-1$, $j=m-2$, and $0\leq j \leq m-3$. We obtain
\begin{align}
  T_g &\leq C \sum\limits_{m=3}^{\infty}\sum\limits_{j=0}^{m-3} {m \choose j} \psi_j \langle D(u u)\rangle_{1,m-j-1} \frac{\tau^{m-3}}{(m-3)!^s}\notag\\
  & + C \sum\limits_{m=3}^{\infty}\Big( \psi_m \Vert u u \Vert_{L^2(\Omega)}+ m \psi_{m-1} \Vert D(u u) \Vert_{L^2(\Omega)}+ m^2 \psi_{m-2}\Vert D^2( u u) \Vert_{L^2(\Omega)}\Big) \frac{\tau^{m-3}}{(m-3)!^s}
\end{align}Using the bound $\psi_j \leq C (j-3)!^s/\tau_{*}^{j}$, the combinatorial estimate $  {m \choose j} (j-3)!^s (m-j-3)!^s/(m-3)!^s \leq C$, and $\eta = \tau/\tau_* <1$, we obtain
\begin{align}
  T_g &\leq C\sum\limits_{m=3}^{\infty} \sum\limits_{j=0}^{m-3} \eta^j \Big( \langle D(uu)\rangle_{1,m-j-1} \frac{\tau^{m-j-3}}{(m-j-3)!^s}\Big) + C \Vert u u \Vert_{H^2(\Omega)}\notag\\
  & \leq C \sum\limits_{m=3}^{\infty} \langle D(uu)\rangle_{1,m-1} \frac{\tau^{m-3}}{(m-3)!^s} + C \Vert uu \Vert_{H^2(\Omega)},
  \end{align}for some sufficiently large constant $C = C(\gamma)$. We now estimate the terms $T_{f}^1$ and $T_{f}^2$. We have
\begin{align*}
  T_{f}^1 &= C \sum\limits_{m=3}^{\infty} \sum\limits_{\alpha_3 = 1}^{m} \epsilon^{\alpha_3} \langle\partial_{3}^{\alpha_3-1} f\rangle_{\alpha_3,m-\alpha_3} \frac{\tau^{m-3}}{(m-3)!^s} \notag\\
   &  \leq C \sum\limits_{m=3}^{\infty} \sum\limits_{|\alpha|=m-1} \epsilon^{\alpha_3+1} \sup\limits_{0<\delta\leq \delta_0} \delta^{m-3} \Vert \partial^\alpha f \Vert_{L^2(\Omega_\delta)} \frac{\tau^{m-3}}{(m-3)!^s},
\end{align*}and similarly
\begin{align*}
  T_{f}^2 &= C \sum\limits_{m=3}^{\infty} \sum\limits_{\alpha_3 = 2}^{m} \epsilon^{\alpha_3} \langle \partial_{3}^{\alpha_3-2} f\rangle_{\alpha_3-1,m-\alpha_3+1} \frac{\tau^{m-3}}{(m-3)!^s} \notag\\
  & \leq C \sum\limits_{m=3}^{\infty} \sum\limits_{|\alpha|=m-1} \epsilon^{\alpha_3+2} \sup\limits_{0<\delta\leq \delta_0} \delta^{m-3} \Vert \partial^\alpha f \Vert_{L^2(\Omega_\delta)} \frac{\tau^{m-3}}{(m-3)!^s}.
\end{align*}Recall that cf.~\eqref{eq:fdef} we have $f = \partial_i u_j \partial_k  u_l D_{ijkl}$, where $D_{ijkl}$ is of Gevrey-class $s$ (cf.~\eqref{eq:Ddef}), and therefore we have $\psi_j \leq C (j-2)!^s/\tau_{*}^{j}$, for all $j\geq 0$. Here we have denoted $\psi_{\beta,\delta} = \delta^{\max\{|\beta|-2,0\}} \Vert \partial^\beta D_{ijkl} \Vert_{L^\infty(\Omega_\delta)}$, and also $\psi_j = \sum_{|\beta|=j} \sup_{0<\delta\leq \delta_0} \psi_{\beta,\delta}$. From the above estimates and the Leibniz rule we obtain that $T_{f}^1 + T_{f}^2$ is bounded by
\begin{align}
  & C \sum\limits_{m=3}^{\infty} \sum\limits_{j=0}^{m-1} {m-1 \choose j} \sum\limits_{|\alpha|=m-1} \sum\limits_{|\beta|=j, \beta\leq \alpha} \epsilon^{\alpha_3 + 1} \sup\limits_{0<\delta\leq \delta_0} \delta^{m-3} \Vert \partial^{\beta} D_{ijkl} \Vert_{L^\infty(\Omega_\delta)} \Vert \partial^{\alpha-\beta}( Du Du) \Vert_{L^2(\Omega_\delta)} \frac{\tau^{m-3}}{(m-3)!^s}\notag\\
  &\qquad \qquad  \leq C \sum\limits_{m=3}^{\infty}\sum\limits_{j=0}^{m-3} {m-1 \choose j} \psi_j \sum\limits_{\alpha_3=0}^{m-j-1}\epsilon^{\alpha_3+1} \langle \partial_3^{\alpha_3}(Du Du)\rangle_{\alpha_3+3,m-j-1-\alpha_3} \frac{\tau^{m-3}}{(m-3)!^s}\notag\\
  &\qquad \qquad \qquad \qquad + C \sum\limits_{m=3}^{\infty}( \psi_{m-1} + m \psi_{m-2} )  \Vert Du Du \Vert_{H^1(\Omega)} \frac{\tau^{m-3}}{(m-3)!^s}.\label{eq:Tf:bound}
\end{align}Using the bound $\psi_j \leq C (j-2)!^s/\tau_{*}^{j}$, the combinatorial estimate
\begin{align}
  {m-1 \choose j} \frac{(j-2)!^s (m-j-3)!^s}{(m-3)!^s} \leq C,
\end{align}and the fact that $\eta = \tau/\tau_{*} < 1$, from \eqref{eq:Tf:bound} we obtain
\begin{align}
  T_{f}^1 + T_{f}^2 & \leq C \sum\limits_{m=3}^{\infty} \sum\limits_{j=0}^{m-3} \eta^j \left(\sum\limits_{\alpha_3=0}^{m-j-1} \epsilon^{\alpha_3+1} \langle\partial_3^{\alpha_3} (Du Du)\rangle_{\alpha_3+3,m-j-1-\alpha_3} \frac{\tau^{m-j-3}}{(m-j-3)!^s}\right) + C \Vert Du Du \Vert_{H^1(\Omega)}\notag\\
  & \leq C \sum\limits_{m=3}^{\infty}\sum\limits_{\alpha_3=0}^{m-1} \epsilon^{\alpha_3+1} \langle \partial_{3}^{\alpha_3} (Du Du)\rangle_{\alpha_3+3,m-1-\alpha_3} \frac{\tau^{m-3}}{(m-3)!^s} + C \Vert Du Du \Vert_{H^1(\Omega)},
\end{align}for some sufficiently large $C = C(\gamma)>0$. This concludes the proof of the lemma. \end{proof}

\subsection*{Proof of Lemma~\ref{lemma:pressureterm}}
Here we use the estimate obtained in Lemma~\ref{lemma:P:u} to bound ${\mathcal P}$ in terms of the Gevrey-class norm of the velocity, and prove the estimate \eqref{eq:pressureterm}. In view of Lemma~\ref{lemma:P:u}, we need to estimate the terms
\begin{align}
  &\sum\limits_{m=3}^{\infty}\sum\limits_{\alpha_3=0}^{m-1} \epsilon^{\alpha_3+1} \langle \partial_3^{\alpha_3} (Du Du)\rangle_{\alpha_3+3,m-1-\alpha_3} \frac{\tau^{m-3}}{(m-3)!^s} \notag\\
  &\qquad \leq C \sum\limits_{m=3}^{\infty} \sum\limits_{j=0}^{m-1}\sum\limits_{|\alpha|=m-1} \sum\limits_{|\beta|=j, \beta\leq \alpha} {\alpha \choose \beta} \epsilon^{\alpha_3 + 1} \sup\limits_{0<\delta\leq \delta_0} \delta^{m-3} \Vert \partial^\beta Du \partial^{\alpha-\beta} Du \Vert_{L^2(\Omega_\delta)} \frac{\tau^{m-3}}{(m-3)!^s} = {\mathcal R}, \label{eq:R}
\end{align}and the lower order term
\begin{align}
  &\sum\limits_{m=3}^{\infty} \langle D(uu)\rangle_{1,m-1} \frac{\tau^{m-3}}{(m-3)!^s}\notag\\
  &\qquad \leq C \sum\limits_{m=3}^{\infty} \sum\limits_{j=0}^{m}\sum\limits_{|\alpha|=m, \alpha_3\leq 1} \sum\limits_{|\beta|=j, \beta\leq \alpha}{\alpha \choose \beta} \sup\limits_{0<\delta\leq \delta_0} \delta^{m-3} \Vert \partial^\beta u \partial^{\alpha-\beta} u \Vert_{L^2(\Omega_\delta)} \frac{\tau^{m-3}}{(m-3)!^s}= {\mathcal S}. \label{eq:S}
\end{align}Similarly to the estimate for the the commutator term ${\mathcal C}$ (cf.~Proof of Lemma~\ref{lemma:velocitycommutator}), bounding ${\mathcal R}$ and ${\mathcal S}$ is achieved by splitting the above sums according to the relative sizes of $j$ and $m-j$. This idea was introduced in our previous work \cite{KV2}.
Namely, we write the right side of \eqref{eq:R} as ${\mathcal R}_1 + {\mathcal R}_2 + {\mathcal R}_3 + {\mathcal R}_{\rm low} + {\mathcal R}_{\rm high} + {\mathcal R}_4 + {\mathcal R}_5$, according to $j=0,1,2$, $3\leq j \leq [(m-1)/2]$, $[(m-1)/2]+1\leq j \leq m-3$, $j=m-2$, and respectively $j=m-1$. Note that by symmetry (replace $j$ by $m-j$) the terms ${\mathcal R}_1$ and ${\mathcal R}_5$, ${\mathcal R}_2$ and ${\mathcal R}_4$, and also ${\mathcal R}_{\rm low}$ and ${\mathcal R}_{\rm high}$, have the same upper bounds. We have the estimates
\begin{align}
{\mathcal R}_1 + {\mathcal R}_5 &\leq C \Vert Du \Vert_{L^\infty(\Omega)} \Vert u \Vert_{H^3(\Omega)} + C\tau \Vert Du \Vert_{L^\infty(\Omega)} \Vert  u \Vert_{Y_\tau} \label{eq:R15}\\
{\mathcal R}_2 + {\mathcal R}_4 &\leq C\Vert D^2 u \Vert_{L^\infty(\Omega)} \Vert u \Vert_{H^2(\Omega)} + C\tau \Vert D^2 u \Vert_{L^\infty(\Omega)} \Vert u \Vert_{H^3(\Omega)} + C\tau^2 \Vert D^2 u \Vert_{L^\infty(\Omega)} \Vert u \Vert_{Y_\tau}\label{eq:R24} \\
{\mathcal R}_3 & \leq C\tau^2 \Vert u \Vert_{H^5(\Omega)} \Vert u \Vert_{Y_\tau},\label{eq:R3}
\end{align}and also
\begin{align}
  {\mathcal R}_{\rm low} + {\mathcal R}_{\rm high} \leq C (\tau^{3/2} + (1+K^3)\tau^{3}) \Vert u \Vert_{X_\tau} \Vert u \Vert_{Y_\tau}.\label{eq:Rmed}
\end{align}The proofs of \eqref{eq:R15}--\eqref{eq:Rmed} are similar to those in \cite[Section 5]{KV2} and those in Section~\ref{sec:commutator} of the present paper, and are thus omitted. Combined they give the desired estimate on ${\mathcal P}$. To estimate ${\mathcal S}$ one proceeds similarly. Note though that this is a lower order term. We have the following bound
\begin{align}
  {\mathcal S} &\leq C \Big( \tau \Vert u \Vert_{L^\infty(\Omega)} + \tau^2 \Vert D u \Vert_{L^\infty(\Omega)} + \tau^3 \Vert  D^2 u \Vert_{L^\infty(\Omega)} + (\tau^{5/2} + (1+K^3)\tau^4) \Vert u \Vert_{X_\tau}\Big) \Vert u \Vert_{Y_\tau}\notag\\
  &\qquad + C(1+\tau^2) \Big( \Vert u \Vert_{W^{2,\infty}(\Omega)}^2 + \Vert u \Vert_{H^3(\Omega)}^2\Big),\label{eq:S}
\end{align}where $C>0$ is a constant that may depend on $\gamma$. The proof of \eqref{eq:S} is omitted (see \cite[Section 5]{KV2} for details). By collecting the above estimates, and the lower order terms from \eqref{eq:P:u}, we conclude the proof of the pressure estimate.

\section{Global Gevrey-class persistence} \label{sec:presistence} \setcounter{equation}{0}
In this section we prove that the local, short time estimates of Section~\ref{sec:short-time-Gevrey} may be combined together to obtain global (in space) Gevrey-class a priori estimates that are valid for all $t<T_*$, the maximal time of existence of the Sobolev solution.

Let $T<T_*$ be fixed. We shall prove that the solution $u(t)$ is of Gevrey-class $s$ on $[0,T]$ and give a lower bound on the radius of Gevrey-class regularity. For this purpose let $\{ x^{\lambda}\}_{\lambda=1}^{N}$ be points on $\partial D$ determined as follows. In a small neighborhood of $x^{\lambda}$ the boundary of $D$ is the graph of a Gevrey-class function $\gamma^{\lambda}$, i.e., there exists $r^{\lambda}>0$ sufficiently small such that $D^\lambda= D \cap B_{r^{\lambda}}(x^{\lambda}) = \{ x\in B_{r^{\lambda}}(x^{\lambda})\, :\, x_3 > \gamma^{\lambda}(x_1,x_2)\}$. Moreover, we can pick $r^{\lambda}$ small enough so that after composing with a rigid body rotation about $x^{\lambda}$ we have $\Vert  \partial_1 \gamma^{\lambda} \Vert_{L^\infty} + \Vert  \partial_2 \gamma^{\lambda} \Vert_{L^\infty} \leq \overline{\varepsilon}$, where $\overline{\varepsilon}>0$ is the fixed universal constant of Lemma~\ref{lemma:H^2}. For all $\lambda\in \{1,\ldots,N\}$ we let $\Omega^{\lambda} = D \cap B_{r^{\lambda}/2}(x^{\lambda})$. We take $N$ large enough so that there exists a compactly embedded open set $\Omega \subset D$ with analytic boundary, such that $\Omega \cup \bigcup_{1\leq \lambda\leq N} \Omega^\lambda = D$. To obtain Gevrey-class regularity in the interior of $D$, we cover $\Omega$ with finitely many, sufficiently small, analytic charts $\{D^\lambda\}_{N+1}^{N+N'}$, chosen as follows. Denote by $\Omega^\lambda$ a ball inside $D^\lambda$, and let $r^\lambda = {\rm dist}(\bar{\Omega^\lambda},(D^\lambda)^c)$, where $\lambda \in \{ N+1,\ldots,N+N'\}$. We let $N'$ be large enough so that
\begin{align}
  1 \leq \sum\limits_{\lambda=1}^{N} \chi_{\Omega^{\lambda}}(x) + \sum\limits_{\lambda=N+1}^{N+N'} \chi_{\Omega^{\lambda}}(x)\leq C \label{eq:covering}
\end{align}for all $x\in D$, where $C\geq 1$ is a sufficiently large constant.

For $s\geq 0$ fixed, define by $\phi_{t,s}(a)$ the particle trajectory with initial condition $\phi_{s,s}(a) = a$, i.e., the unique smooth solution to
\begin{align*}
  &\frac{d}{dt} X(t) = u(X(t),t)\\
  & X(s) = a.
\end{align*}Note that $\phi_{t,0}(a) = \phi_t(a)$, where $\phi_t$ is as defined in \eqref{eq:X1}--\eqref{eq:X2}. Since the flow map $\phi_{t,s}\colon D\mapsto D$ is a bijection, cf.~\eqref{eq:covering}, we also have $1\leq \sum_{\lambda = 1}^{N+N'} \chi_{\phi_{t,s}(\Omega^\lambda)}(x) \leq C$ for all $0\leq s \leq t$ and all $x\in D$.

Let $T_0=0$, and define $T_1$ as the maximal time $0=T_0 < T_1 \leq T$ such that for all $T_0\leq t \leq T_1$ we have that $\phi_{t,T_0}(\Omega^\lambda) \subset D^\lambda$ for all $\lambda \in \{1,\ldots,N+N'\}$. Note that if $T_1< T$, then by the maximality of $T_1$, there exists $\lambda \in \{1,\ldots,N+N'\}$ with $\phi_{T_1,T_0}(\Omega^{\lambda}) \cap \overline{(D^{\lambda})^{c}} \neq \emptyset$. Thus there exists and $x_{0}\in \Omega^\lambda$ such that $| \phi_{T_1,T_0}(x_{0}) - x_{0}| \geq r^\lambda/2 \geq r^*$, where $r^* = \min_{1\leq \lambda \leq N+N'} \{r^\lambda/2\}$ is a fixed constant. We obtain that if $T_1 < T$, then $T_1$ may be estimated from below via
\begin{align}
  \int_{T_{0}}^{T_1} \Vert u(\cdot, t) \Vert_{W^{1,\infty}(D)}\, dt \geq r^*.\label{eq:T1estimate}
\end{align} For each $\lambda \in \{1,\ldots,N+N'\}$, let $\theta^\lambda (x_1,x_2,x_3) = (x_1,x_2,x_3 - \gamma^\lambda(x_1,x_2))=(y_1,y_2,y_3)$ be a boundary straightening map and define $\tilde{\Omega}^\lambda = \theta^\lambda(\Omega^\lambda)$. Note that this is exactly the setup from Section~\ref{sec:notation}. Let $u^\lambda(x,t) = u(x,t) \chi_{D^\lambda}(x)$ and for $y = \theta^\lambda(x) \in \theta^\lambda(D^\lambda)$ define $\tilde{u}^\lambda(y,t) = u^\lambda(x,t)$.

Let $\tau_0=\tau(T_0)$ be the uniform radius of Gevrey-class regularity of the initial data $u_0$. By possibly decreasing $\tau_0$ by a factor, we may assume that $\tau_0\leq \epsilon \tau_*$, where $\epsilon=\epsilon(D)>0$ is as in Lemma~\ref{lemma:mathcal:P}, and $\tau_*$ is the uniform  radius of Gevrey-class regularity of $\partial D$. Since $u^\lambda(T_0)$ has Gevrey-class radius $\tau_0$, we have that $\Vert \tilde{u}^\lambda(T_0,y) \Vert_{X_{a_* \tau_0}} < \infty$ for all $\lambda \in \{1,\ldots,N+N'\}$, where $0<a_*\leq 1$ measures the possible decrease in the Gevrey-class radius after composing with the boundary straightening map $\theta^\lambda$ (cf.~Remark~\ref{rem:FaadiBruno}). Therefore, on $[T_0,T_1]$ we can apply Theorem~\ref{thm:short:time} for each chart $\{\Omega^\lambda\}_{\lambda=1}^{N}$, respectively Remark~\ref{rem:interior} for $\{\Omega^\lambda\}_{\lambda=N+1}^{N+N'}$, to obtain that for all $\lambda\in\{1,\ldots,N+N'\}$ we have (cf.~\eqref{eq:Gevrey:growth})
\begin{align}
\Vert \tilde{u}^\lambda(\cdot,t) \Vert_{X_{\tau(t)}} \leq Q_0 + C \int_{T_0}^{t} \Big(1+K^2(s)\Big) M^2(s)\; ds \label{eq:Q1estimate}
\end{align}where $Q_0 = \max_{\lambda\in\{1,\ldots,N+N'\}} \Vert \tilde{u}^\lambda(\cdot,T_0) \Vert_{X_{a_* \tau_0}}$, $C = C(D)$ is a positive constant, and the radius of Gevrey-class regularity $\tau(t)$ is bounded from below (cf.~\eqref{eq:tau:compact}) by
\begin{align}
  \tau(t) \geq a_* \tau_0\Big( 1 + C t  Q_0 + C t^2 M^2(T_0) \Big)^{-2} \exp\Big(C K(T_0) - C K(t)\Big)\label{eq:it0:tau}
\end{align}for all $T_0 \leq t \leq T_1$. Here we recall that $K(t)=\int_{0}^{t} \Vert u(\cdot,s) \Vert_{W^{1,\infty}(D)}\; ds$, and $M(t) = \Vert u(\cdot,t) \Vert_{H^r(D)}$. Therefore, modulo composing with $(\theta^\lambda)^{-1}$ we obtain that the localized velocity $u^\lambda(x,t)$ is of Gevrey-class $s$ on $[T_0,T_1]$ for each $\lambda \in \{1,\ldots,N+N'\}$. By \eqref{eq:covering} we obtain that $u(\cdot,t)$ is of Gevrey-class $s$ on $[T_0,T_1]$, with uniform radius of Gevrey-class regularity bounded from below by $a_*$ times the right side of \eqref{eq:it0:tau}.

We proceed inductively. Let $k \geq 1$ be fixed. Since $\phi_{T_k,T_{k-1}}(D)=D$ , as above for $t=0$ we cover $D$ with local charts $\{\Omega^\lambda\}_{\lambda=1}^{N+N'}$ and define $T_{k+1}$ as the maximal time $T_{k+1}\leq T$ such that $\phi_{t,T_{k}}(\Omega^\lambda) \subset D_\lambda$ for all $\lambda\in \{1,\ldots,N+N'\}$. Similarly to \eqref{eq:T1estimate} we obtain that if $T_{k+1} < T$, then $T_{k+1}$ may be estimated from
\begin{align}\label{eq:Tk}
\int_{T_k}^{T_{k+1}} \Vert u(\cdot, t) \Vert_{W^{1,\infty}(D)}\, dt \geq r^*.
\end{align}

The induction assumption is that $u(x,T_k)$ if of Gevrey-class $s$, the uniform (over $x\in D$) radius of Gevrey-class regularity of $u(x,T_k)$ is bounded from below by
\begin{align}\label{eq:Tkestimate}
\tau_k = a_{*}^{2} &\; \tau_{k-1} \Big( 1 + C (T_k-T_{k-1})  Q_{k-1} + C (T_{k}-T_{k-1})^{2} M^2(T_{k-1}) \Big)^{-2} \exp\Big(C K(T_{k-1}) - C K(T_k)\Big),
\end{align}and that the Gevrey-class norm at $t=T_k$, given by $Q_k = \max_{\lambda\in\{1,\ldots,N+N'\}} \Vert \tilde{u}^\lambda(\cdot,T_k) \Vert_{X_{a_* \tau_k}}$, is bounded as
\begin{align}
   Q_k \leq Q_{k-1} + C \int_{T_{k-1}}^{T_k} \Big(1+K^2(s)\Big) M^2(s)\; ds. \label{eq:Qkestimate}
\end{align}
We apply Theorem~\ref{thm:short:time}, respectively Remark~\ref{rem:interior}, on each local chart $\Omega^\lambda$, and conclude that $\tilde{u}^{\lambda}(y,t)$ is of Gevrey-class $s$ on $\in[T_k,T_{k+1}]$ for all $\lambda\in\{1,\ldots,N+N'\}$, with Gevrey-class norm bounded as
\begin{align}
  \Vert \tilde{u}^\lambda(\cdot,t) \Vert_{X_{\tau(t)}} \leq Q_k + C \int_{T_k}^{t} \Big(1+K^2(s)\Big) M^2(s)\; ds \label{eq:Qk+1estimate},
\end{align}and radius of Gevrey-class regularity $\tau(t)$ bounded from below by
\begin{align}
  a_{*}\; \tau_k \Big( 1 + C (t-T_k)  Q_{k} + C (t-T_k)^{2} M^2(T_{k}) \Big)^{-2} \exp\Big(C K(T_{k}) - C K(t)\Big). \label{eq:Tk+1estimate}
\end{align}Modulo composing with the inverse map of $\theta^\lambda$, if follows from \eqref{eq:covering} and \eqref{eq:Qk+1estimate} that $u(x,t)$ is of Gevrey-class $s$ for all $t\in[T_k,T_{k+1}]$ with radius bounded from below by $a_*$ times the quantity in \eqref{eq:Tk+1estimate}. Moreover, letting $t=T_{k+1}$ in \eqref{eq:Qk+1estimate}--\eqref{eq:Tk+1estimate} we obtain that the induction assumptions \eqref{eq:Tkestimate}--\eqref{eq:Qkestimate} hold for the next iteration step.

We claim that for each fixed $T<T_*$ the inductive argument described above stops after finitely many steps, i.e., there exists a $k\geq 1$ such that $T_k = T$. To see this, note that if $T_k < T$, then from \eqref{eq:T1estimate} and \eqref{eq:Tk} we obtain
\begin{align}
k r_* \leq \int_{0}^{T_k} \Vert u(\cdot,t) \Vert_{W^{1,\infty}(D)}\; dt \leq \int_{0}^{T} \Vert u(\cdot,t) \Vert_{W^{1,\infty}(D)}\; dt < \infty,\label{eq:k:estimate}
\end{align}which cannot hold for all $k\geq 1$, proving the claim. Moreover, we proved that it takes at most $[K(T)/r_*] +1$ applications of Theorem~\ref{sec:short-time-Gevrey} to show that $u(\cdot,T)$ is uniformly of Gevrey-class $s$, where $[\cdot]$ denotes the integer part, and $K(t)$ is as usual defined by \eqref{eq:K:def}.

It is left to prove that the uniform radius of Gevrey-class regularity $\tau(T)$ of $u(\cdot,T)$ depends explicitly on the initial data and $K(T)$. Let $k = [K(T)/r_*] +1$ and hence $T=T_k$. It follows form the above paragraph that $\tau(T) \geq \tau_{k}$. By the induction assumptions \eqref{eq:Tkestimate}--\eqref{eq:Qkestimate} we bound $\tau_k$ from below as
\begin{align}
  \tau_k &\geq a_{*}^{2k} \tau_0 \prod\limits_{j=1}^{k} \exp\Big(C K(T_{j-1}) - C K(T_j)\Big)  \Big( 1 + C (T_j-T_{j-1})  Q_{j-1} + C (T_{j}-T_{j-1})^{2} M^2(T_{j-1}) \Big)^{-2}.
\end{align}Since $a_{*}^{2k} \leq \exp( - 2k \log(1/a_*) ) \leq \exp(- 2 K(T) \log(1/a_*)/r_*)$ we obtain that
\begin{align}
  \tau_k & \geq \tau_0 \exp\Big( - C K(T)\Big)  \prod\limits_{j=1}^{k} \Big( 1 + C (T_{j}-T_{j-1})  Q_{j-1} + C (T_{j}-T_{j-1})^{2} M^2(T_{j-1}) \Big)^{-2}
\end{align}for a sufficiently large constant $C$ depending only on the domain. To estimate the product term in the above inequality we note that by \eqref{eq:Qkestimate} we have that $Q_{j-1} \leq Q_0 + C M^2(0) \exp(C K(T_{j-1}))$, while from the Sobolev energy estimate we obtain $M^2(T_{j-1}) \leq M^2(0) \exp(C K(T_{j-1}))$. Therefore we have
\begin{align*}
  \tau_k & \geq \tau_0 \exp\Big( - C K(T)\Big)  \prod\limits_{j=1}^{k} \left( 1 + C (T_{j}-T_{j-1})  Q_{0} + C (T_{j}-T_{j-1}) (1+T) M^2(0) \exp\Big(C K(T_{j-1})\Big)  \right)^{-2}\notag\\
  & \geq \tau_0 \exp\Big( - C K(T)\Big) \exp \Big( - C \sum\limits_{j=1}^{k} K(T_{j-1})\Big) \prod\limits_{j=1}^{k} \Big( 1 + C (T_{j}-T_{j-1})  Q_{0} + C (T_{j}-T_{j-1}) (1+T) M^2(0) \Big)^{-2}.
\end{align*}By using the inequality between the arithmetic and the geometric mean, and the fact that $k = C K(T)$, we obtain
\begin{align}
  \tau(T) & \geq \tau_0 \exp\Big( -C \sum\limits_{j=1}^{k} K(T_j) \Big) \left(1 + \frac{C T Q_0 + C T^2 M^2(0)}{k}\right)^{-2k}\notag\\
  & \geq \tau_0 \exp\Big( -C K^2(T) \Big) \exp \Big( - C T Q_0 - C T^2 M^2(0) \Big).
\end{align}
Therefore we have proven the following statement, which is the main theorem of this paper.
\begin{theorem}\label{thm:main}
  Let $u_0$ be divergence-free and of Gevrey-class $s$, with $s\geq 1$, on a Gevrey-class $s$, open, bounded domain $D \subset {\mathbb R}^3$, and $r\geq 5$. Then the unique solution $u(\cdot,t) \in C([0,T_*);H^r(D))$ to the initial value problem \eqref{eq:E1}--\eqref{eq:E4} is of Gevrey-class $s$ for all $t<T_*$, where $T_* \in (0,\infty]$ is the maximal time of existence in $H^r(D)$. Moreover, the radius $\tau(t)$ of Gevrey-class regularity of the solution $u(\cdot,t)$ satisfies
\begin{align}\label{eq:thm:main}
  \tau(t)\geq  C \tau_0  \exp\left(- C \left(\int_{0}^{t} \Vert u(s) \Vert_{W^{1,\infty}} ds\right)^2 \right) \exp\Big(- C t \Vert u_0 \Vert_{X_{\tau_0}} - C t^2 \Vert u_0 \Vert_{H^r}^2\Big),
\end{align}for all $t<T_*$, where $C$ is a sufficiently large constant depending only on the domain $D$, $\tau_0$ is the radius of Gevrey-class regularity of the initial data $u_0$, and $\Vert u_0 \Vert_{X_{\tau_0}}$ is its Gevrey-class norm.
\end{theorem}
\begin{remark}
  Theorem~\ref{thm:main} also holds in the case of a two-dimensional Gevrey-class domain. In 2D it is known that $\Vert u(s) \Vert_{W^{1,\infty}} \leq C \exp(C t)$ for some positive constant $C = C(D,u_0)$, and therefore estimate \eqref{eq:thm:main} shows that the radius of Gevrey-class regularity of the solution is bounded from below by $C \exp(-C \exp(Ct))$ for some $C>0$, depending on the domain and on the initial data. We note that such a lower bound on $\tau(t)$ was obtained in the $2D$ analytic case $s=1$ by Bardos, Benachour, and Zerner \cite{BBZ}, whereas in the non-analytic Gevrey-class case on domains with generic boundary, Theorem~\ref{thm:main} is the first such result (see also \cite{KV1} for the periodic domain, and \cite{KV2} for the half-plane).
\end{remark}

\appendix
\section{}\label{sec:appendix}\setcounter{equation}{0}
\begin{lemma}\label{lemma:combinatorial}
Let $\{a_\lambda \}$, and $\{b_{\lambda,\mu} \}$ be sequences of positive numbers, where $\lambda,\mu \in {\mathbb N}_0^3$. The identity
\begin{align*}
  \sum\limits_{|\alpha|=m} \sum\limits_{|\beta|=j,\; \beta \leq \alpha} \sum\limits_{|\gamma|=k,\; \gamma \leq \beta} a_\gamma b_{\alpha-\beta,\beta-\gamma} = \left(\sum\limits_{|\gamma|=k} a_\gamma \right) \left(\sum\limits_{|\alpha|=m-k} \sum\limits_{|\beta|=j-k,\; \beta\leq \alpha} b_{\alpha-\beta,\beta}\right)
\end{align*}holds for positive integers $j,k,m$ such that $k\leq j \leq m$.
\end{lemma}
\begin{proof}
  \begin{align}
    \sum\limits_{|\alpha|=m} \sum\limits_{|\beta|=j,\; \beta \leq \alpha} \sum\limits_{|\gamma|=k,\; \gamma \leq \beta} a_\gamma b_{\alpha-\beta,\beta-\gamma} &= \sum\limits_{|\alpha|=m}\sum\limits_{|\gamma|=k,\; \gamma \leq \alpha}a_\gamma \sum\limits_{|\beta|=j,\; \gamma \leq \beta \leq \alpha} b_{(\alpha -\gamma)-(\beta-\gamma), \beta-\gamma}\notag \\
    &=\sum\limits_{|\alpha|=m}\sum\limits_{|\gamma|=k,\; \gamma \leq \alpha}a_\gamma \sum\limits_{|\lambda|=j-k,\; \lambda \leq \alpha-\gamma} b_{(\alpha-\gamma)-\lambda,\lambda}  \notag\\
    & = \sum\limits_{|\alpha|=m}\sum\limits_{|\gamma|=k,\; \gamma \leq \alpha}a_\gamma d_{\alpha-\gamma}\label{eq:ap:lemma:1},
  \end{align}where we let
  \begin{align*}
    d_{\alpha-\gamma} = \sum\limits_{|\lambda|=j-k,\; \lambda \leq \alpha-\gamma} b_{(\alpha-\gamma)-\lambda,\lambda}.
  \end{align*}By \cite[Lemma 4.2]{KV2} the far right side of \eqref{eq:ap:lemma:1} may be written as
  \begin{align*}
    \left(\sum\limits_{|\gamma|=k} a_\gamma \right) \left( \sum\limits_{|\alpha|=m-k} d_\alpha\right) = \left(\sum\limits_{|\gamma|=k} a_\gamma \right) \left( \sum\limits_{|\alpha|=m-k} \sum\limits_{|\beta|=j-k,\; \beta\leq \alpha} b_{\alpha-\beta, \beta }\right),
  \end{align*}which concludes the proof of the lemma.
\end{proof}

\begin{lemma}\label{lemma:ap:convolve}
Let $0<\eta < 1$ and $\{a_{m,j}\}_{m\geq0,j\geq0}$ be a sequence of positive numbers. Then we have
\begin{align}\label{eq:ap:conv}
\sum\limits_{m=3}^{\infty} \sum\limits_{j=1}^{m} \sum\limits_{k=0}^{j} \eta^k a_{m-k,j-k} &= \frac{\eta^3}{1-\eta} a_{0,0} + \frac{\eta^2}{1-\eta} (a_{1,0} + a_{1,1}) + \frac{\eta}{1-\eta} (a_{2,0}+a_{2,1}+ a_{2,2}) \notag\\
&\qquad  + \frac{\eta}{1-\eta} \sum\limits_{m=3}^{\infty} a_{m,0} + \frac{1}{1-\eta} \sum\limits_{m=3}^{\infty}\sum\limits_{j=1}^{m} a_{m,j}.
\end{align}
\end{lemma}
\begin{proof}
By re-indexing we have
\begin{align*}
  \sum\limits_{m=3}^{\infty} \sum\limits_{j=1}^{m} \sum\limits_{k=1}^{j} \eta^k a_{m-k,j-k} &= \sum\limits_{m=3}^{\infty} \sum\limits_{k=1}^{m} \eta^k \left(\sum\limits_{j=k}^{m} a_{m-k,j-k}\right) =  \sum\limits_{m=3}^{\infty} \sum\limits_{k=1}^{m} \eta^k \left(\sum\limits_{j=0}^{m-k} a_{m-k,j}\right) = \sum\limits_{m=3}^{\infty} \sum\limits_{k=1}^{m} \eta^k b_{m-k},
\end{align*}where we denoted $b_l = \sum_{j=0}^{l} = a_{l,j}$. By summing the geometric series in $\eta$ the far right side of the above equality may be re-written as $( \eta^3 b_0 + \eta^2 b_1 + \eta \sum_{j=2}^{\infty} b_l)/(1-\eta)$. Therefore we obtain
\begin{align*}
  \sum\limits_{m=3}^{\infty} \sum\limits_{j=1}^{m} \sum\limits_{k=0}^{j} \eta^k a_{m-k,j-k} &= \sum\limits_{m=3}^{\infty} \sum\limits_{j=1}^{m} a_{m,k} + \frac{\eta^3}{1-\eta} a_{0,0} + \frac{\eta^2}{1-\eta} (a_{1,0} + a_{1,1}) + \frac{\eta}{1-\eta} \sum\limits_{m=2}^{\infty} \sum\limits_{j=0}^{m} a_{m,j},
\end{align*}and \eqref{eq:ap:conv} follows by grouping appropriate terms.
\end{proof}

\begin{lemma}\label{lemma:ap:1}
  Let $\{F_\delta(t)\}_{\delta\in[0,\delta_0]}$ be a family of nonnegative $C^1$ functions, where $\delta_0 \leq 1$ is a fixed constant. Assume that
   \begin{enumerate}
   \item $\{\dot{F}_\delta(t)\}_{\delta\in[0,\delta_0]}$ is a uniformly equicontinuous family,
    \item for every fixed $t$, the functions $F_\delta(t)$ and $\dot{F}_\delta(t)$ depend continuously on $\delta$.
    \end{enumerate}Then for every fixed $t\in(0,\infty)$ we have
  \begin{align}
   \frac{d^+}{dt} \sup\limits_{\delta\in [0,\delta_0]} F_\delta(t) =  \limsup\limits_{h\rightarrow 0^+} \frac 1h \left( \sup\limits_{\delta\in[0,\delta_0]} F_\delta(t+h) - \sup\limits_{\delta\in[0,\delta_0]} F_\delta(t)\right) \leq \sup\limits_{\delta\in[0,\delta_0]} \dot{F}_\delta(t).\label{eq:1}
  \end{align}
\end{lemma}
\begin{proof}
  Fix $t\in (0,\infty)$. For a fixed $\delta \in [0,\delta_0]$, and $h>0$, we have
  \begin{align*}
    F_\delta(t+h) = \int_{t}^{t+h} \dot{F}_\delta(s)\, ds + F_\delta(t) \leq \int_{t}^{t+h} \left(\sup\limits_{\delta\in[0,\delta_0]} \dot{F}_\delta(s)\right)\, ds + \sup\limits_{\delta\in[0,\delta_0]} F_\delta(t).
  \end{align*}Therefore
  \begin{align*}
  \frac 1h \left(\sup\limits_{\delta\in[0,\delta_0]} F_\delta(t+h) - \sup\limits_{\delta\in[0,\delta_0]}F_\delta(t)\right)
 \leq \frac{1}{h} \int_{t}^{t+h} \left(\sup\limits_{\delta\in[0,\delta_0]} \dot{F}_\delta(s)\right)\, ds,  \end{align*}and if we can prove that $\sup_{\delta\in[0,\delta_0]} \dot{F}_\delta(t)$ is a continuous function of $t$, then \eqref{eq:1} holds, concluding the proof of the lemma. The fact that $\sup_{\delta\in[0,\delta_0]} \dot{F}_\delta(t)$ is a continuous function of $t$ follows directly from the definition of uniform equicontinuity and the inequality
 \begin{align*}
   \left| \sup\limits_{\delta\in[0,\delta_0]} a(\delta) - \sup\limits_{\delta\in[0,\delta_0]} b(\delta) \right| \leq    \sup\limits_{\delta\in[0,\delta_0]} | a(\delta) - b(\delta) |,
 \end{align*}which holds for all functions $a,b:[0,\delta_0] \rightarrow {\mathbb R}$.
\end{proof}

\begin{lemma}\label{lemma:ap:2}
  Let $\tilde{v} = \partial^\alpha \tilde{u}$, for some $\alpha \in {\mathbb N}_{0}^{3}$, and $\tilde{u}$ as in Lemma~\ref{lemma:energyestimate}. Let
  \begin{align}\label{eq:ap:2}
    f_\delta(t) = \delta^{|\alpha|} \Vert \tilde{v}(t,\cdot) \Vert_{L^2(\tilde{\Omega}_{\delta,t})},
  \end{align}and $F_\delta(t) = f_{\delta}^{2}(t)$. Then the family $\{F_{\delta}(t)\}_{\delta\in[0,\delta_0]}$ satisfies the conditions {\rm (i)} and {\rm (ii)} of Lemma~\ref{lemma:ap:1}.
\end{lemma}
\begin{proof}
  Let $\tilde{\Omega}_t=\bigcup_{\delta\in[0,\delta_0]} \tilde{\Omega}_{\delta,t}$, and $\phi_t(x)$ be the particle trajectory with initial data $x$. Without loss of generality assume that $\tilde{\Omega}_t \in \tilde{D}$ for all $t>0$, and that $\Vert \tilde{v}\Vert_{L^2(\tilde{\Omega})} \neq 0$.

  The fact that for a fixed $t$ the family $F_{\delta}(t)$ depends continuously on $\delta$, follows from the continuity of the integral with respect to the Lebesgue measure, and the fact that $\tilde{v} \in L^\infty(\tilde{D})$. Also, from \eqref{eq:E2} and the fact that $\det(\partial \theta / \partial x) = 1$, we have
  \begin{align}\label{eq:3}
    \frac{d}{dt} \Vert \tilde{v}(t,\cdot) \Vert_{L^2(\tilde{\Omega}_{\delta,t})}^2 = 2 \int\limits_{\Omega_\delta} g(t,\phi_t(x)) v(t,\phi_t(x))\, dx,
  \end{align}and since $gv \in L^\infty(D)$, the continuity of the integral implies that $\dot{F}_\delta(t)$ depends continuously on $\delta$, so that condition (ii) holds.

  To show that the family $\{ \dot{F}_\delta(t)\}$ is uniformly equicontinuous, let $\epsilon>0$. We need to show that there exists $\tau=\tau(\epsilon)>0$ such that $|\dot{F}_\delta(t) - \dot{F}_\delta(s)| < \epsilon$ for all $|t-s|\leq \tau$ and all $\delta\in[0,\delta_0]$. By \eqref{eq:3} and the mean value theorem we have
  \begin{align*}
  &\left| \frac{d}{dt} \Vert \tilde{v}(t,\cdot) \Vert_{L^2(\tilde{\Omega}_{\delta,t})}^2 - \frac{d}{dt} \Vert \tilde{v}(s,\cdot) \Vert_{L^2(\tilde{\Omega}_{\delta,s})}^2 \right|\\
    & \qquad \qquad = 2 \left| \int\limits_{\Omega_\delta} \Big( g(t,\phi_t(x)) v(t,\phi_t(x)) - g(s,\phi_s(x)) v(s,\phi_s(x)) \Big)\, dx \right|\\
    & \qquad \qquad \leq C |t-s| \sup\limits_{z\in(t,s)} \sup\limits_{x\in\Omega} \left| \partial_t (gv) (z,\phi_z(x)) \right| + \left| u_j(z,\phi_z(x)) \partial_j(gv)(z,\phi_z(x))\right|\\
    & \qquad \qquad \leq C |t-s| \left( \Vert \partial_t (gv) \Vert_{L^\infty(D)}  + \Vert u \Vert_{L^\infty(D)} \Vert gv  \Vert_{W^{1,\infty}(D)}\right),
  \end{align*}since $\tilde{\Omega}_t \in \tilde{D}$ for all $t>0$. Recall that $\tilde g = [\partial^\alpha, \tilde{u}_j \, \partial_j \theta_k\, \partial_k] \tilde{u} + \partial^\alpha (\partial_j \theta_k\, \partial_k \tilde{p})$, and so the right side of the above is bounded by $C |t-s| \Vert u \Vert_{H^{r+|\alpha|}(D)}^3$ for some sufficiently large $r$. To conclude the proof of the lemma one follows standard arguments.\end{proof}

  \begin{lemma}\label{lemma:ap:3}
    Let $\tilde{v} = \partial^\alpha \tilde{u}$, for some $\alpha \in {\mathbb N}_{0}^{3}$, and $\tilde{u}$ as in Lemma~\ref{lemma:energyestimate}. Let $M_\delta(t) \geq 0$ be an upper bound
    \begin{align}
      \delta^{|\alpha|} \frac{d}{dt} \Vert \tilde{v}(t,\cdot) \Vert_{L^2(\tilde{\Omega}_{\delta,t})} \leq M_\delta(t)\label{eq:ap:A5}
    \end{align}which holds for all $t\geq 0$ and all $\delta\in[0,\delta_0]$ such that $\Vert \tilde{v}(t,\cdot) \Vert_{L^2(\tilde{\Omega}_{\delta,t})} > 0$. Furthermore, assume that $\sup_{\delta\in[0,\delta_0]} M_\delta(t)$ is continuous in $t$. Then we have
    \begin{align}
      \frac{d^+}{dt\,} \sup\limits_{\delta\in[0,\delta_0]} \delta^{|\alpha|} \Vert \tilde{v}(t,\cdot) \Vert_{L^2(\tilde{\Omega}_{\delta,t})} \leq \sup\limits_{\delta\in[0,\delta_0]} M_\delta(t),
    \end{align}where we denote by $d^+a(t)/dt  = \limsup_{h\rightarrow 0+}( a(t+h)-a(t) )/h$ the right derivative of a function.
  \end{lemma}
  \begin{proof}
    Let $f_\delta=\delta^{|\alpha|} \Vert \tilde{v}(t,\cdot) \Vert_{L^2(\tilde{\Omega}_{\delta,t})}$ and let $F_\delta = f^2_\delta$. Note that by assumption we have $\dot{f}_{\delta}(t) \leq M_\delta(t)$ for all $\delta \in [0,\delta_0]$ and $t>0$. It follows from Lemmas~\ref{lemma:ap:1} and~\ref{lemma:ap:2} that
    \begin{align*}
      \frac{d^+}{dt\,} \sup\limits_{\delta\in[0,\delta_0]} F_\delta(t) \leq \sup\limits_{\delta\in[0,\delta_0]} \dot{F}_\delta(t).
    \end{align*}Due to the continuity in $\delta$ of the family $f_\delta$, $$\sup_{\delta\in[0,\delta_0]} F_\delta = \left(\sup_{\delta\in[0,\delta_0]} f_\delta\right)^2.$$ Therefore,
    \begin{align}
      \frac{d^+}{dt\,} \sup\limits_{\delta\in[0,\delta_0]} f_\delta(t) = \frac{\frac{d^+}{dt\,} \sup\limits_{\delta\in[0,\delta_0]} F_\delta(t)}{2 \sup\limits_{\delta\in[0,\delta_0]} f_\delta(t)} \leq \frac{\sup\limits_{\delta\in[0,\delta_0]} f_\delta(t) M_\delta(t)}{\sup\limits_{\delta\in[0,\delta_0]} f_\delta(t)} \leq \sup_{\delta\in[0,\delta_0]} M_\delta(t),\label{ap:eq:sup=0}
    \end{align}for all $t$ such that $\sup_{\delta\in[0,\delta_0]} f_\delta(t) > 0$. This concludes the proof of the lemma by noting that if $g(t)$ is a nonnegative function such that $(d^+/dt) g(t) \leq G(t)$ for all $t$ such that $g(t)>0$, with $G(t)$ continuous, then $g(t) \leq g(t_0) + \int_{t_0}^{t} G(s)\; ds$, for all $0\leq t_0 < t$.  \end{proof}

\begin{ack}
  The work of both authors was supported in part by the NSF grant DMS-1009769.
\end{ack}

\end{document}